\documentclass[a4paper, 11pt]{article} 
\usepackage{amsmath, amscd, amsthm, amssymb, mathrsfs}

\DeclareMathOperator{\im}{im}
\DeclareMathOperator{\re}{re}

\DeclareMathOperator{\End}{End}

\DeclareMathOperator{\tr}{tr}

\DeclareMathOperator{\diag}{diag}

\DeclareMathOperator{\SO}{SO}
\DeclareMathOperator{\U}{U}
\DeclareMathOperator{\SU}{SU}
\DeclareMathOperator{\GL}{GL}
\DeclareMathOperator{\SL}{SL}
\DeclareMathOperator{\PSL}{PSL}

\DeclareMathOperator{\Rm}{Rm}
\DeclareMathOperator{\Ric}{Ric}
\DeclareMathOperator{\ric}{ric}

\DeclareMathOperator{\dvol}{dvol}

\newcommand{\R}{\mathbb R}
\newcommand{\C}{\mathbb C}
\newcommand{\Z}{\mathbb Z}

\newcommand{\diff}{\text{\rm d}}
\newcommand{\del}{\partial}

\newcommand{\so}{\mathfrak{so}}
\newcommand{\su}{\mathfrak{su}}
\newcommand{\gl}{\mathfrak{gl}}

\newcommand{\g}{\mathfrak{g}}

\renewcommand{\P}{\mathbb P}

\theoremstyle{plain}
	\newtheorem{theorem}{Theorem}[section]
	\newtheorem{proposition}[theorem]{Proposition}
	\newtheorem{lemma}[theorem]{Lemma}
	\newtheorem{corollary}[theorem]{Corollary}

\theoremstyle{definition}
	\newtheorem{definition}[theorem]{Definition}
	\newtheorem{remark}[theorem]{Remark}

\theoremstyle{plain}
	\newtheorem*{theorem*}{Theorem}
	\newtheorem*{proposition*}{Proposition}
	\newtheorem*{lemma*}{Lemma}
	\newtheorem*{corollary*}{Corollary}
	\newtheorem*{conjecture*}{Conjecture}
\theoremstyle{definition}
	\newtheorem*{definition*}{Definition}
	\newtheorem*{remark*}{Remark}
	\newtheorem*{remarks*}{Remarks}
		
\numberwithin{equation}{section}

\begin{document}

\title{Symplectic Calabi--Yau manifolds, minimal surfaces and the hyperbolic geometry of the conifold}
\author{
Joel Fine\footnote{Supported by an FNRS postdoctoral fellowship.} \ and Dmitri Panov\footnote{Supported by EPSRC grant EP/E044859/1.}
}

\maketitle

\vfill
\begin{abstract}
Given an $\SO(3)$-bundle with connection, the associated two-sphere bundle carries a natural closed 2-form. Asking that this be symplectic gives a curvature inequality first considered by Reznikov \cite{reznikov}. We study this inequality in the case when the base has dimension four, with three main aims. 

Firstly, we use this approach to construct symplectic six-manifolds with $c_1=0$ which are never K\"ahler; e.g., we produce such manifolds with $b_1=0=b_3$ and also with $c_2\cdot [\omega] <0$, answering questions posed by Smith--Thomas--Yau \cite{smith-thomas-yau}. 

Examples come from Riemannian geometry, via the Levi--Civita connection on $\Lambda^+$. The underlying six-manifold is then the twistor space and often the symplectic structure tames the Eells--Salamon twistor almost complex structure. Our second aim is to exploit this to deduce new results about minimal surfaces: if a certain curvature inequality holds, it follows that the space of minimal surfaces (with fixed topological invariants) is compactifiable; the minimal surfaces must also satisfy an adjunction inequality, unifying and generalising results of Chen--Tian \cite{chen-tian}. 

One metric satisfying the curvature inequality is hyperbolic four-space $H^4$. Our final aim is to show that the corresponding symplectic manifold is symplectomorphic to the small resolution of the conifold $xw-yz=0$ in $\C^4$. We explain how this fits into a hyperbolic description of the conifold transition, with isometries of $H^4$ acting symplectomorphically on the resolution and isometries of $H^3$ acting biholomorphically on the smoothing. 
\end{abstract}
\vfill
\newpage

\vfill
\setcounter{tocdepth}{2}
\tableofcontents
~
\vfill

\newpage

\section{Introduction}

\subsection{Symplectic Calabi--Yau six-manifolds}

The first main aim of this article is to describe a construction of symplectic non-K\"ahler Calabi--Yau six-manifolds. For the purposes of this article, by ``symplectic Calabi--Yau'' we mean the following:

\begin{definition}
A \emph{symplectic Calabi-Yau manifold} is a symplectic manifold with $c_1=0$.
\end{definition}

(Note that simple connectivity is \emph{not} part of this definition; indeed we don't present any simply connected examples here, although some do have $b_1=0$, see \S\ref{fake projective planes}.)

There are, of course, already known examples of symplectic Calabi--Yau manifolds which are not K\"ahler, the most famous, perhaps, being the Kodaira--Thurston torus \cite{thurston}. This example is a Nil-manifold (a quotient of a nilpotent Lie group by a cocompact lattice) and this idea has been exploited by several authors \cite{gray,gray2,benson-gordon} to produce further examples of symplectic non-K\"ahler Calabi--Yau manifolds. In dimension four, there is a folk-lore conjecture that essentially all minimal non-K\"ahler symplectic Calabi--Yaus arise this way. Some evidence for this is provided by the recent work of Li \cite{li} who obtains strong restrictions on the rational cohomology of a minimal symplectic four-manifold with $c_1$ torsion. 

In dimension six, the nilpotent Lie groups admitting left-invariant symplectic forms have been classified by Salamon \cite{salamon}. He finds twenty-six such groups all of which can be quotiented by a cocompact lattice to give twenty-six types of Nil-manifold which support symplectic Calabi--Yau structures, none of which are K\"ahler. By comparison, the method for constructing symplectic Calabi--Yau manifolds described here seems closer in spirit to hyperbolic geometry. The examples given show that in dimension six there are more symplectic Calabi--Yaus than are accounted for by Nil-manifolds.

\subsection{Definite connections in dimension four}

The construction central to this article produces symplectic forms on certain sphere bundles. Given an $\SO(3)$-bundle $E \to X$ with compatible connection $\nabla$, the associated two-sphere bundle $Z \to X$ carries a natural closed 2-form $\omega$. Asking that $\omega$ be symplectic gives a curvature inequality for $\nabla$ which was first studied by Reznikov \cite{reznikov}. (The inequality is explained in detail in \S\ref{basic construction}.)

Here we focus on the curvature inequality---and the corresponding symplectic manifolds---in the case when the base $X$ has dimension four. In this dimension, the inequality has a simple geometric formulation: the curvature of $\nabla$ on any two-plane in $TX$ must be non-zero. We call connections satisfying this condition \emph{definite}. 

Another feature of the four-dimensional situation is that definite connections carry a sign. The distinction between positive and negative definite connections is fundamental: positive definite connections yield symplectic Fano manifolds (i.e., ones for which $[\omega]$ is a positive multiple of $c_1$) whilst negative definite connections yield symplectic Calabi--Yau manifolds. (These cohomological computations are carried out in \S\ref{cohomology}.)

Unfortunately, we have been unable to use this construction to produce any new symplectic Fano manifolds. There are, however, large numbers of negative definite connections. The main source of examples comes from Riemannian geometry: on an oriented Riemannian four-manifold, the Levi--Civita connection induces an $\SO(3)$-connection on the bundle of self-dual 2-forms $\Lambda^+$ and it makes sense to ask for it to be definite. This gives an inequality for the Riemann tensor. (When the metric is anti-self-dual, the inequality yielding positive definite connections was first studied by Gauduchon \cite{gauduchon}; we describe the relationship between his approach and that considered here in \S\ref{anti-self-dual solutions}.)

As is explained in \S\ref{connections from metrics}, there are many solutions to the inequality which give negative definite connections; examples include hyperbolic and complex-hyperbolic metrics, metrics with pointwise $2/5$-pinched negative sectional curvatures and certain metrics constructed by Gromov--Thurston \cite{gromov-thurston} with pointwise $K$-pinched negative sectional curvatures with $K$ arbitrarily large.  Whilst this list may suggest that the inequality is intimately linked with negative sectional curvature, there are differences between the two situations. For example, there are complete solutions to the inequality on the total space of $\mathcal O(-n)\to \C\P^1$ when $n\geq3$. This manifold has $\pi_2\neq 0$ and so supports no complete metrics of negative sectional curvature.  On the other hand, we do not know of an example of a definite connection on a \emph{compact} manifold with non-trivial higher homotopy; in particular, all compact examples we have found so far have infinite fundamental group. We do not yet know if this is a necessary condition. 

The Calabi--Yau manifolds produced in this way are not K\"ahler.
Hodge theoretic considerations show that no compact symplectic Calabi--Yau manifold coming from a negative definite connection can have a compatible complex structure. These manifolds don't just violate K\"ahlerian Hodge theory, however. Using this construction we produce examples of symplectic Calabi--Yau six-manifolds with $b_1=0=b_3$ and with $c_2\cdot[\omega]<0$ (see \S\ref{asd solutions}); both of these topological situations are impossible for a K\"ahler Calabi--Yau. These examples answer two questions posed by Smith--Thomas--Yau \cite{smith-thomas-yau}.

\subsection{Twistors and minimal surfaces in four-manifolds}

The second main aim of this article is to explore the relationship between twistor geometry and the symplectic geometry of definite connections arising from Riemannian metrics; this is done in \S\ref{twistor spaces}. When the Levi--Civita connection on $\Lambda^+$ is definite, the unit sphere bundle $Z$ of $\Lambda^+$ carries a symplectic form $\omega$. The manifold $Z$---the twistor space of $X$---also carries two natural almost complex structures: the Atiyah--Hitchin--Singer almost complex structure $J_+$ \cite{atiyah-hitchin-singer} and the Eells--Salamon almost complex structure $J_-$ \cite{eells-salamon}. Asking for $\omega$ to tame $J_+$ or $J_-$ gives two further curvature inequalities which are special cases of those ensuring the connection on $\Lambda^+$ is positive or negative definite respectively. Again, in the negative case, this inequality has many solutions; except for the Gromov--Thurston examples, for all the metrics listed above, $\omega$ tames $J_-$.

The almost complex structure $J_-$ has special significance for minimal surfaces: a result of Eells--Salamon \cite{eells-salamon} gives a one-to-one correspondence between minimal surfaces (i.e., branched minimal immersions of surfaces) in $X$ and non-vertical pseudoholomorphic curves in $(Z, J_-)$. When $\omega$ tames $J_-$ we can use this observation to deduce facts about minimal surfaces. For example, Gromov compactness implies  that certain spaces of minimal surfaces have natural compactifications (see \S\ref{gromov compactness}). It is also possible to use the symplectic form to prove an adjunction inequality: when $\omega$ tames $J_-$, any embedded minimal surface $\Sigma \subset X$ must satisfy
$$
\chi(\Sigma) + \Sigma \cdot \Sigma < 0.
$$
This is proved in \S\ref{adjunction} (which also covers branched minimal immersions). In the case of minimal surfaces in hyperbolic manifolds this inequality goes back to Wang. For the more general metrics considered here this result unifies and generalises some of the adjunction inequalities of Chen--Tian \cite{chen-tian}.

To give these results some context, \S\ref{twistor spaces} begins with a discussion of some classical results concerning minimal surfaces in \emph{three}-dimensional manifolds. When the ambient space $M$ is negatively curved, minimal immersions $\Sigma \to M$ must have $\chi(\Sigma)<0$ and, moreover, the space of branched minimal immersions of $\Sigma$ is compact. These two facts can also be seen from a symplectic point of view and from this perspective their direct analogues are the four-dimensional results presented here. One interesting point is that the four-dimensional inequality is genuinely different from negative curvature: there are curvature tensors which satisfy the four-dimensional inequality which have some sectional curvatures positive.

\subsection{Hyperbolic geometry and the conifold}

The third main aim of this article, addressed in \S\ref{quadric cone},  is to show how hyperbolic geometry features when one considers the threefold quadric cone $Q=\{xw-yz\}\subset \C^4$---also known as the conifold---and its desingularisations. As is mentioned above, there is a natural symplectic structure on the twistor space of $H^4$. We show in \S\ref{symplecto to small res} that the twistor space is symplectomorphic to the total space of $R = \mathcal O(-1)\oplus \mathcal O (-1)\to \C\P^1$. (The symplectic structure on $R$ is obtained by adding the pull-backs of the standard symplectic forms from $\C\P^1$ and $\C^4=\C^2\oplus\C^2$ via the obvious maps.)

The manifold $R$ arises as one way to remove the singularity in the conifold $Q$; the map $R \to \C^4$ has image $Q$ and exhibits $R$ as the so-called small resolution. An alternative choice of desingularisation of $Q$ is provided by the smoothing $S = \{xw-yz = 1\}$. These two different desingularisations have been much studied in both the mathematics and physics literature (for example, \cite{friedman,lu-tian,smith-thomas-yau,strominger,candelas-delaossa}). The passage from one desingularisation to the other is called a ``conifold transition'' and, loosely speaking, the two sides $R$ and $S$ of the conifold transition are considered to be a ``mirror pair'' (see various articles in \cite{essays}, or the paper \cite{morrison}). 

The description of $R$ as the twistor space of $H^4$ shows that isometries of $H^4$ act by symplectomorphisms on $R$. On the other side of the picture, the group of biholomorphisms of $S$ is $\SO(4,\C)$. Since $\SO(4,\C) \cong \SL(2,\C)\times \SL(2,\C)/\pm1$ and $\PSL(2,\C)$ is the group of isometries of $H^3$, we see that \emph{three}-dimensional hyperbolic geometry features strongly in the complex geometry of $S$. The fact that $R$ and $S$ are mirror to each other makes one naively wonder if there is any hidden relationship between three- and four-dimensional hyperbolic geometry. In \S\ref{mirror symmetry} we describe certain superficial similarities between the complex geometry of quotients of $S$ and the symplectic geometry of quotients of $R$, although this picture still remains disappointingly vague.

\subsection{Acknowledgments}

During the course of this work we have benefited from conversations with many people. It is a pleasure to acknowledge the comments and advice of Martin Bridson, Tom Coates, Simon Donaldson, Mark Haskins, Maxim Kontsevich, Anton Petrunin, Simon Salamon, Michael Singer, Ivan Smith, Richard Thomas and Henry Wilton. We would also like to thank Gideon Maschler for bringing the article of Reznikov \cite{reznikov} to our attention. 

\section{Symplectic forms from connections}

\subsection{The basic construction}\label{basic construction}

Let $E \to X^{2n}$ be an $\SO(3)$-bundle over a $2n$-manifold. Write $\pi\colon Z \to X$ for the associated unit sphere bundle and $V \to Z$ for the vertical tangent bundle of $Z$. A choice of metric connection $\nabla$ in $E$ induces both a splitting $TZ = V \oplus H$ and a metric connection $\nabla^V$ in $V$ defined as follows: in the vertical directions $\nabla^V$ is the Levi--Civita connection of the fibres; in the horizontal directions, parallel transport with respect to $\nabla$ isometrically identifies nearby fibres and hence connects their tangent spaces.

Since $\nabla^V$ is a metric connection, its curvature $2\pi i\omega$ is an imaginary two-form. (Strictly speaking, to write the curvature in this way we need also to choose an orientation on each fibre of $Z$, changing the orientation simply swaps the sign of $\omega$.) The restriction of $\omega$ to a fibre is, of course, just the area form of the fibre. The horizontal components of $\omega$, however, encode the curvature $F$ of the original connection $\nabla$ in a certain way.

To describe this, denote by $h \colon \so(3) \to C^\infty(S^2)$ the map which associates to each $u \in \so(3)$ its corresponding mean-value zero Hamiltonian function $h(u)$. Explicitly, $2\pi h(u)(p) = (u, u_p)$, where $u_p$ denotes the unit-speed right-handed rotation about the axis in the direction of $p \in S^2$ and $(\cdot , \cdot)$ denotes the standard inner-product in $\so(3)$.

Using this ``comoment map'' fibrewise gives a map, $h\colon\Omega^0_X(\so(E)) \to C^\infty(Z)$.  Tensoring this with pull-back $\pi^*\colon \Omega^p_X \to \Omega^p_Z$ extends $h$ to a map on endomorphism valued forms $h \colon \Omega^p_X(\so(E)) \to \Omega^p_Z$. With this definition in hand, the following Proposition gives the exact relationship between $\omega$ and $F$. 

\begin{proposition}\label{decomposition of omega}
With respect to the decomposition $\Lambda^2_Z \cong \Lambda^2 V^* \oplus (V^* \otimes H^*) \oplus \Lambda^2H^*$,
$$
\omega = \omega_{S^2}\oplus 0\oplus h(F),
$$
where $2\pi\omega_{S^2}$ is the fibrewise area form. In particular, $\omega$ is a symplectic form precisely when $h(F)^n$ is nowhere zero.
\end{proposition}

\begin{proof}
By construction, the $\Lambda^2V^*$-component of $\omega$ is, up to a factor of $2\pi$, the area form of each fibre. The $V^*\otimes H^*$ component measures the failure of $\nabla$-parallel transport to commute with Levi--Civita parallel transport in the fibre. Since, however, $\nabla$-parallel transport is an isometry on the fibres of $Z$, this component is zero. It remains to compute the $\Lambda^2H^*$-component.

Let $h,k$ be vector fields on $X$, let $\tilde h,\tilde k$ denote their horizontal lifts to $Z$ and let $v$ be a vertical field on $Z$; we wish to compute 
$$
F_{\nabla^V}(\tilde h, \tilde k)(v)
=
\left(
\nabla^V_{\tilde k} \nabla^V_{\tilde h} 
-
\nabla^V_{\tilde h} \nabla^V_{\tilde k}
+
\nabla^V_{[\tilde h,\tilde k]}
\right)v.
$$
The horizontal component of $[\tilde h, \tilde k]$ is the horizontal lift of $[h,k]$, which we denote $\widetilde{[h,k]}$, whilst the vertical component is $u = F(h, k)$ interpreted as a vertical vector field. Since vertically $\nabla^V$ is the Levi--Civita connection $\nabla^{S^2}$ of the fibres we have
\begin{equation}\label{curvature equation}
F_{\nabla^V}(\tilde h,\tilde k)(v)
=
\left(
\nabla^V_{\tilde k} \nabla^V_{\tilde h} 
-
\nabla^V_{\tilde h} \nabla^V_{\tilde k}
+
\nabla^V_{\widetilde{[h,k]}}
\right) v
+
\nabla^{S^2}_u v.
\end{equation}

Now, the original connection $\nabla$ induces a connection $\nabla'$ on the vector-bundle over $X$ whose fibre at $x\in X$ is the space $C^\infty(T \pi^{-1}(x))$ of all vector fields on the 2-sphere over at $x$. The bracketed terms in (\ref{curvature equation}) are precisely the curvature $F_{\nabla'}(h,k)(v)$ of $\nabla'$; this is given by interpreting the curvature $u=F(h,k)$ of $\nabla$ as acting on $C^\infty(T\pi^{-1}(x))$:
$$
\left(
\nabla^V_{\tilde k} \nabla^V_{\tilde h} 
-
\nabla^V_{\tilde h} \nabla^V_{\tilde k}
+
\nabla^V_{\widetilde{[h,k]}}
\right) v
=
[v,u].
$$
(This is an application of the standard theory of connections in principal bundles: such a connection induces connections in all associated bundles; given a representation $\rho \colon G \to \GL(W)$, the curvature of the induced connection in the associated $W$-bundle is given by composing the curvature of the original connection---which takes values in $\g$---with the map $
\rho_*\colon \g \to \gl(W)$.)

Hence
$
F_{\nabla^V}(h,k)(v)
=
[v,u]
+
\nabla^{S^2}_u v
$.
Since $\nabla^{S^2}$ is torsion-free, this simplifies to $F_{\nabla^V}(h,k)(v) = \nabla^{S^2}_v u$. This reduces our curvature calculation to the following question purely about the geometry of $S^2$: given $u \in \so(3)$, the map $T_pS^2 \to T_pS^2$ given by $v \mapsto \nabla^{S^2}_v u$ is a rotation by what angle? The following lemma ensures that this is $2\pi h(u)(p)$ where $h(u)$ is the Hamiltonian of $u$, giving the claimed decomposition of $\omega$.

As it is a curvature form, $\omega$ is automatically closed, hence symplectic if and only if $\omega^{n+1} = \omega_{S^2} \wedge h(F)^n$ is nowhere vanishing. This happens precisely when $h(F)^n$ is nowhere zero.
\end{proof}

\begin{lemma}
Let $u \in \so(3)$. The map $T_pS^2 \to T_pS^2$ given by $v \mapsto \nabla^{S^2}_v u$ is a rotation by $2\pi h(u)(p)$.
\end{lemma}

\begin{proof}
Think of $u$ as a vector field on $\R^3$; it is given simply by multiplication: $x\mapsto u(x)$.  Since $u(x)$ is linear in $x$, $\nabla^{\R^3}_v u = u(v)$. The connection on $S^2$ is induced by projection, giving $\nabla^{S^2}_v u =  u(v) - (u(v), p)\, p$. In other words, $v \mapsto \nabla_v^{S^2}u$ is given by the component of $u$ which is rotation about the axis through $p$. The size of this component is exactly $2\pi h(u)$.
\end{proof}

\begin{remark}\label{high dim fibres}~
\begin{enumerate}
\item
An alternative way to compute the $\Lambda^2H^*$-component of $\omega$ is to use a result appearing, for example, in Chapter 6 of \cite{mcduff-salamon}. This result states that if $\Omega$ is a closed 2-form on a fibre bundle $M\to B$ whose fibrewise restriction is non-degenerate, it defines a horizontal distribution $H$;  given vectors $h,k$ at $p\in B$ with horizontal lifts $\tilde h, \tilde k$ the function $\Omega(\tilde h, \tilde k)$ defined on the fibre $M_p$ over $p$ is a Hamiltonian for the curvature $F_H(h,k)$ with respect to the symplectic form $\Omega|_{M_p}$. This fixes $\Omega(\tilde h, \tilde k)$ up to the addition of a constant. In other words, the $\Lambda^2H^*$-component of $\Omega$ is determined up to the addition of a form pulled-back from $B$. Applying this in our case we see that $\omega$ is precisely as claimed, at least modulo $\Lambda^2 T^*X$. 
\item
There is a similar result for bundles with fibres other than $S^2$. Since this article focuses on two-sphere bundles we only give a brief description. Let $G$ be a Lie group which admits a symplectic quotient $(G/H, \omega_{G/H})$. Suppose, moreover, that this symplectic form is the curvature of a unitary connection on a Hermitian line bundle $L\to G/H$ and that the action of $G$ on $G/H$ lifts to a unitary action of $G$ on $L$ preserving the connection. (This is equivalent to the existence of a moment map $\mu \colon G/H \to \g^*$, embedding $G/H$ as an integral coadjoint orbit.) 

Now, given a principal $G$-bundle $P\to X$ with connection $\nabla$, let $Z \to X$ denote the associated bundle with fibres $G/H$ and $\mathcal L \to Z$ the Hermitian line bundle whose fibrewise restriction is $L\to G/H$. The connection $\nabla$ enables one to combine the fibrewise connections in $\mathcal L$ to give a unitary connection $\nabla^{\mathcal L}$. The curvature of this connection is
$2\pi i\,\omega$ where, again using the decomposition $\Lambda^2V^* \oplus (V^*\otimes H^*) \oplus \Lambda^2 H^*$ defined by $\nabla$,
$$
\omega
=
\omega_{G/H} \oplus 0 \oplus h(F)
$$
in which $F$ is the curvature of $\nabla$ and $h(F)$ is defined as before, using the map $h\colon \g \to C^\infty(G/H)$ dual to the moment map $\mu$, i.e., with $h(\xi)(q) = \langle \mu(q), \xi\rangle$.
\end{enumerate}
\end{remark}

\subsection{Connections over four-manifolds}

When $X$ is an orientable four-manifold, the condition that $\omega$ be symplectic has a simple geometric formulation. Consider the curvature $F$, a section of  $\Lambda^2T^*X \otimes \so(E)$, as a map $F \colon \so(E)^* \to \Lambda^2T^*X$. In dimension four, the wedge product $\Lambda^2\times \Lambda^2 \to \Lambda^4$ defines an indefinite conformal structure on $\Lambda^2$ of signature $(3,3)$. Using this, the condition $h(F)^2 \neq 0$ has an alternative interpretation: $h(F)^2$ is nowhere vanishing if and only if the image of $F$ is a maximal definite subspace of $\Lambda^2$. 

Alternatively, consider the dual map $F^* \colon \Lambda^2TX \to \so(E)$ which associates to each bivector $u\wedge v$ the holonomy of $A$ in the plane spanned by $u$ and $v$. Saying that $F$ has  maximal definite image in $\Lambda^2T^*X$ is equivalent to saying that the kernel of $F^*$ has trivial intersection with the cone of decomposable vectors in $\Lambda^2TX$. In other words, $\omega$ is symplectic if and only if the connection $\nabla$ has non-zero holonomy on every two-plane in $TX$.

\begin{definition}
An $\SO(3)$-connection over a four-manifold is called \emph{definite} if the image of its curvature map $F \colon \so(E)^* \to \Lambda^2$ is a maximal definite subspace for the wedge product. Equivalently, a connection is definite if and only if it has non-zero holonomy on every  two-plane in $TX$.
\end{definition}

\begin{remark}
Although we don't use it here, we also give, for completeness, the generalisation of this definition to higher dimensional bases and fibres. In the context of the second part of Remark \ref{high dim fibres}, a connection in a principal $G$-bundle $P\to X$ is called \emph{$G/H$-definite} if at every point $x\in X$ the composition of the moment and curvature maps
$$
h(F)(x) \colon (G/H)_x \stackrel{\mu}{\to} \g^*_x \stackrel{F}{\to} \Lambda^2 T^*_xX
$$
has image which misses the ``null cone'' of 2-forms whose top power vanishes. Given a $G/H$-definite connection, the associated $G/H$-bundle is naturally symplectic.
\end{remark}

Returning to $\SO(3)$-connections over four-manifolds, note that a definite connection not only defines a symplectic structure on $Z$ but also an orientation and metric on $X$. The Pontrjagin form $\tr F^2$ is a positive multiple of $\pi_*\omega^3$ (this follows from the description of $\omega$ in Proposition \ref{decomposition of omega}) and so is nowhere vanishing and serves as the volume form. Next, note that given any maximal definite subspace $\Lambda^+\subset \Lambda^2$, there is a unique conformal structure for which $\Lambda^+$ is the bundle of self-dual two forms. Taking $\Lambda^+$ to be the image of $F$ determines a conformal class and hence, together with the volume form, a metric on $X$. Notice this is the unique conformal structure which makes $\nabla$ an instanton; conversely, any instanton whose curvature map has no kernel is a definite connection.

\subsubsection{The sign of a definite connection}

A key feature of definite connections on four-manifolds is that they carry a sign; this is defined as follows. The Lie algebra $\so(3)$ has a natural orientation. (This is analogous to the fact that $\C$ is naturally oriented; there is a choice involved, but once the convention is agreed upon, all copies of $\so(3)$, $\so(3)$-bundles etc.\ are oriented.) One way to see this is to chose an orientation on $\R^3$; the orientation determines an isomorphism $f\colon\R^3 \to \so(3)$ given by taking $f(v)$ to generate right-handed rotation of speed $|v|$ about the directed axis through $v$. Changing the orientation of $\R^3$ changes the sign of $f$ so that the push-forward of the orientation of $\R^3$ by $f$ is independent of the orientation on $\R^3$ that you started with. 

For a definite connection, the curvature map is a bundle isomorphism $\so(E)^* \cong \Lambda^+$. As we've just described, $\so(E)^*$ is an oriented bundle and the same is true for $\Lambda^+$. For example, the metric defined by the definite connection gives an identification $\Lambda^2=\so(TX)$ and $\Lambda^+\subset \Lambda^2$ corresponds to an $\so(3)$ summand in the splitting $\so(4)=\so(3)\oplus\so(3)$. This gives $\Lambda^+$ an orientation. As $F$ is an isomorphism between oriented bundles, we can ask whether $\det F$ is positive or negative.

\begin{definition}
A definite connection is called \emph{positive} if its curvature map $F \colon \so(E)^* \to \Lambda^+$ is orientation preserving, and \emph{negative} if it is orientation reversing.
\end{definition}

This means that the space of definite curvature tensors has four components given by considering the induced orientation and sign. As we will see in the next section, the difference between positive and negative is fundamental; it corresponds to $(Z, \omega)$ being ``Fano'' (i.e., $[\omega]$ a positive multiple of $c_1(Z)$) or ``Calabi--Yau'' (i.e., $c_1(Z) = 0$).

\subsection{K\"ahlerity and cohomological considerations}\label{cohomology}

Let $k = \frac{1}{2}c_1(V) \in H^2(Z, \R)$. It follows from the Leray--Hirsch theorem that $H^*(Z, \R)$ is a free $H^*(X,\R)$-module generated by $k$. The ring structure is determined by the single additional relation  $k^2= \frac{1}{4}p_1(E)$ (pulled back to $Z$). In the case when $E$ admits a definite connection, the curvature map gives an isomorphism $E \cong \Lambda^+$. In particular, this implies that $k^3 = \frac{1}{4}\int_X p_1(\Lambda^+) = \frac{1}{4}(2\chi + 3\tau)$ where $\chi$ is the Euler characteristic of $X$ and $\tau$ the signature (using the orientation on $X$ induced by $\pi_*\omega^3$). Since the symplectic class is $[\omega] = 2k$, this number must be positive. This interpretation of the symplectic volume in terms of the topology of $X$ gives a first obstruction to the existence of a definite connection:

\begin{proposition}\label{half Hitchin-Thorpe}
If $X$ admits a definite connection then $2\chi + 3\tau >0$ (where the signature $\tau$ is defined using the induced orientation).
\end{proposition}

This inequality is ``one half'' of the Hitchin--Thorpe inequality \cite{hitchin,thorpe}. In the situation considered in \S\ref{connections from metrics}, where the definite connection is induced by the Levi--Civita connection of a Riemannian metric, there is a proof of this inequality along more standard lines; see the discussion following Proposition \ref{Riemannian Hitchin-Thorpe}.
  
The concrete description of $H^*(Z, \R)$ also makes it straightforward to check the Hard Lefschetz property. 

\begin{lemma}\label{Hard Lefschetz}
The map $k^2 \colon H^1(Z, \R) \to H^5(Z,\R)$ is zero. In particular, if $H^1(X, \R) \neq 0$ and $X$ admits a definite connection, the symplectic manifold $(Z, \omega)$ admits no compatible K\"ahler structure.
\end{lemma}

\begin{proof}
This is immediate from the fact that $k^2$ is the  pull-back of a top degree class from $X$ whilst $H^1(Z,\R)$ is also pulled back from $X$.
\end{proof}

More topological information is provided by the Chern classes of $(Z, \omega)$. To compute these the first step is to make a choice of compatible almost complex structure. Whilst for a general symplectic manifold there is no natural choice, in the case of a definite connection on a four-manifold, there is one compatible almost complex structure which is singled out. 

\begin{definition}
On the vertical tangent spaces $V$, take $J$ to be the standard complex structure of the sphere. At a point $p\in Z$, the horizontal component of $\omega$ spans a line in $\Lambda^+_{\pi(p)}$. The plane orthogonal to this line is the real projection of $\Lambda^{2,0}$ for a unique almost complex structure on $T_{\pi(p)}X$. We use this almost complex structure to define $J$ on $H_p\cong T_{\pi(p)}X$. Equivalently, the metric determined by the definite connection gives an identification between rays in $\Lambda^+$ and almost complex structures on $X$ compatible with the metric and orientation. Under this identification, the horizontal component of $\omega$ at $p$ defines an almost complex structure on $T_{\pi(p)}X$.
\end{definition}

\begin{remark}
Of course, this definition is motivated by that of the almost complex structure on twistor space. Recall that the twistor space of a four-manifold with a conformal structure is the sphere bundle associated to $\Lambda^+$. A definite connection determines a metric on $X$ and also gives an identification $\so(E)^* \to \Lambda^+$. By identifying $E \cong \so(E) \cong \so(E)^*$, the curvature map gives a fibrewise linear embedding $Z \to \Lambda^+$. Projecting to the sphere bundle gives a diffeomorphism between $Z$ and the twistor space of $X$. It preserves or reverses orientations on the fibres according to whether the connection is positive or negative definite.

The twistor space carries two natural almost complex structures (one for each choice of fibrewise orientation), but there are no direct relationships between these and $\omega$ or $J$. In particular, they are not in general compatible with $\omega$. We will see later, however, that when the definite connection is defined via a Riemannian metric on $X$ it \emph{is} possible to relate the twistor geometry and symplectic geometry. This is discussed in \S\ref{twistor spaces} where the relationship leads to results concerning minimal surfaces.
\end{remark}

Returning to the calculation of the  Chern classes of $(Z, \omega)$, we use the compatible almost complex structure $J$ and look for the Chern classes of the complex vector bundle $(TZ, J)$. The calculation is similar to Hitchin's for twistor spaces \cite{hitchin2}, the only difference coming from a sign change in the negative definite case.

\begin{lemma}
If $\nabla$ is a positive definite connection, the total Chern class of $Z$ is
$$
c(Z)
=
1 + 4 k + (e + 4k^2) + 2k\cdot e,
$$
where $e$ is the Euler class of $X$ pulled back to $Z$. If $\nabla$ is a negative definite connection, the total Chern class of $Z$ is
$$
c(Z)
=
1 + 0 + (e - 4k^2) + 2k\cdot e.
$$
In particular, positive definite connections give symplectic Fano manifolds whilst negative definite connections give symplectic Calabi--Yau manifolds.
\end{lemma}

\begin{proof}
Since $TZ = V \oplus H$, $c(Z) = c(V)\cdot c(H)$. Now $c(V) = 1 + 2k$ and it remains to work out the Chern classes of $H$. 

To compute $c_1(H)$ we use the following claim: if $\nabla$ is positive definite, there is an isomorphism of complex line bundles $\Lambda^2H  \cong V$; meanwhile, if $\nabla$ is negative definite, $\Lambda^2H \cong V^*$. To prove this claim, note that by definition of $J$, the bundle $\Lambda^2H$ is isomorphic to the sub-bundle of $\pi^*\Lambda^+$ given by the orthogonal complements to the horizontal components of $\omega$. This bundle is, in turn, isomorphic to the tangent bundle of the unit spheres in $\Lambda^+$. As explained in the remark concerning twistor spaces above, this is isomorphic to $V$ or $V^*$ according as $\nabla$ is positive or negative definite. Hence $c_1(H) =\pm 2k$ with the sign agreeing with that of $\nabla$.

Next, note that $c_2(H) = e(H)$ and, since $H \cong \pi^* (TX)$ as orientable bundles, this gives $c_2(H) = e$, the Euler class of $X$ pulled back to $Z$. So the total Chern class of $H$ is
$c(H)= 1 \pm 2k + e$, yielding
$$
c(Z) 
=
(1 + 2k)(1\pm 2k + e),
$$
which gives the result.
\end{proof}

We also record for later use some numerical invariants of $(Z, \omega)$:
\begin{corollary}\label{Chern numbers}
A positive definite connection gives $(Z, \omega)$ with
$$
c_1^3 = 16(2\chi + 3\tau),\quad
c_1\cdot c_2 = 12(\chi + \tau), \quad 
c_3 = 2\chi,\quad
c_2 \cdot [\omega] = 6(\chi + \tau).
$$ 
A negative definite connection gives $(Z, \omega)$ with
$$
c_1^3 = 0,\quad
c_1\cdot c_2 = 0, \quad 
c_3 = 2\chi,\quad
c_2 \cdot [\omega] = -2(\chi + 3\tau).
$$ 
\end{corollary}

We saw above (Lemma \ref{Hard Lefschetz}) that when $H^1(X, \R)\neq 0$, the symplectic manifolds arising from definite connections are not K\"ahler.  The following well-known lemma also rules out K\"ahler examples arising from \emph{negative} definite connections when $H^1(X, \R)=0$.

\begin{lemma}\label{holomorphic volume form}
Let $Z$ be a K\"ahler threefold with $b_1=0$ and $c_1=0$. Then $b_3 >0$.
\end{lemma}
\begin{proof}
Since $b_1=0$, it follows from Hodge theory that $h^{0,1}=0$ and so, by the Dolbeault theorem, that $H^1(Z, \mathcal O) = 0$. Considering the long exact sequence in cohomology associated to the short exact sequence
$$
0 \to \Z \to \mathcal O 
\stackrel{\exp}{\longrightarrow} \mathcal{O}^* \to 1
$$
gives that the first Chern class is injective $H^1(Z, \mathcal O^*) \to H^2(Z, \Z)$. Since $c_1=0$ this implies that the canonical bundle is holomorphically trivial. Now a non-zero holomorphic volume form gives a nontrivial element of $H^3(Z, \C)$.
\end{proof}

\begin{corollary}
The symplectic manifolds arising from negative definite connections on compact four-manifolds never admit compatible K\"ahler structures.
\end{corollary}
\begin{proof}
If $H^1(X, \R) \neq 0$, this follows from the failure of Hard Lefschetz (Lemma \ref{Hard Lefschetz}). If $H^1(X, \R)  =0$ then $b_1(Z)= 0 = b_3(Z)$ and so the result follows from the previous Lemma.
\end{proof}

(We will later see that \emph{positive} definite connections can lead to K\"ahler manifolds.) In many examples of negative definite connections, there are other properties beside Hard Lefschetz which fail. We will indicate these where appropriate.

\section{Symplectic forms from Riemannian metrics}
\label{connections from metrics}

\subsection{A Riemannian curvature inequality}

It is possible to use Riemannian geometry to produce many examples of definite connections, a technique first exploited by Reznikov \cite{reznikov}. On an oriented Riemannian four-manifold, the Hodge star gives a map $\ast \colon \Lambda^2\to\Lambda^2$ with $\ast^2=1$. Accordingly, its eigenvalues are $\pm1$ and the bundle of two-forms splits $\Lambda^2=\Lambda^+\oplus \Lambda^-$ into eigenspaces; the eigenvectors are called self-dual and anti-self-dual forms respectively. The Levi--Civita connection induces a metric connection on, say, the $\SO(3)$-bundle $\Lambda^+$ and it makes sense to ask for this to be a definite connection. When this happens, we obtain a symplectic form on the twistor space of $X$. The case when the metric is anti-self-dual---so that the twistor space is a complex threefold---was first considered by Gauduchon using complex geometry \cite{gauduchon}. We describe his work in \S\ref{anti-self-dual solutions}.

Asking for the Levi--Civita connection on $\Lambda^+$ to be definite translates into a Riemannian curvature inequality. To describe it, we begin by recalling the decomposition of the Riemannian curvature of a four-manifold. The curvature can be thought of as a self-adjoint operator $\Rm \colon\Lambda^2 \to \Lambda^2$ and so, with respect to the decomposition $\Lambda^2=\Lambda^+\oplus\Lambda^-$, it decomposes into parts:
$$
\Rm
=
\left(\begin{array}{cc}
W^+ + \frac{s}{12} & \Ric_0^*\\
&\\
\Ric_0 & W^- + \frac{s}{12}
\end{array}\right).
$$ 
Here $W^+$ and $W^-$ are the self-dual and anti-self-dual Weyl curvatures respectively---they are trace-free self-adjoint operators $W^\pm\colon\Lambda^\pm \to \Lambda^\pm$---whilst $s$ is the scalar curvature, acting by multiplication. $\Ric_0$ is the trace-free Ricci curvature interpreted as an operator $\Lambda^+ \to \Lambda^-$. More explicitly, denote by $\ric_0\colon T^*X \to T^*X$ the trace free Ricci curvature acting as an endomorphism of the cotangent bundle; then $\Ric_0\colon \Lambda^2\to \Lambda^2$ is given on decomposable vectors by
\begin{equation}\label{trace-free Ricci operator}
\Ric_0(a\wedge b) 
= 
\frac{1}{2}(\ric_0(a)\wedge b +a \wedge \ric_0(b))
\end{equation}
and this is readily seen to be self-adjoint and swap self-dual and anti-self-dual forms. (In fact, there are other conventions implicit in this description of $\Rm$. For example, the standard isomorphism $\Lambda^2 \cong \so(TX)$ is only intrinsic up to sign; we use the map $\Lambda^2\otimes TX \subset T^*X\otimes (T^*X \otimes TX) \to T^*X$ given by contracting on the bracketed $T^*X\otimes TX$ factor, followed by the metric isomorphism $T^*X\to TX$ to explicitly define the map $\Lambda^2 \to \so(TX)$. Another convention needed explicitly for later calculation is that fixing the scale of the inner-product on $\Lambda^2$; we declare that given an orthonormal basis $e_i$ of $T^*X$, the bivectors $e_i\wedge e_i$ for $i<j$ are an orthonormal basis of $\Lambda^2$.)

To describe the curvature of the Levi--Civita connection on $\Lambda^+$, use the metric and natural orientation to identify $\Lambda^+ \cong \so(\Lambda^+)^*$, so that the curvature of $\Lambda^+$ is a map $\Lambda^+ \to \Lambda^2$. With this identification, the curvature is simply given by the first column of the decomposition of the Riemann tensor above. In other words, 
$$
F=
\left(W^+ + \frac{s}{12}\right) \oplus \Ric_0 
\colon 
\Lambda^+ \to \Lambda^+\oplus \Lambda^-.
$$
With this in hand we can characterise those metrics giving definite connections on $\Lambda^+$:

\begin{theorem}\label{Riemannian curvature inequality}
Let $X$ be an oriented Riemannian four-manifold. The Levi--Civita  connection on $\Lambda^+$ is definite if and only if the endomorphism of $\Lambda^+$ given by  
\begin{equation}\label{main operator}
\mathcal D
=
\left(W^++\frac{s}{12}\right)^2
-
\Ric_0^*\Ric_0
\end{equation}
is definite. 

We split into two cases according to the induced orientation:
\begin{enumerate}
\item
If $\mathcal D$ is positive definite then the connection induces the original orientation on $X$. In this case the sign of the connection agrees with the sign of $\det(W^++\frac{s}{12})$.

\item
If $\mathcal D$ is negative definite then the connection induces the opposite orientation on $X$. In this case the sign of the connection agrees with the sign of $\det \Ric_0$.
\end{enumerate}

\end{theorem}

\begin{proof}
The connection is definite if the wedge product is definite on the image of $F$. For $v \in \Lambda^+$, 
\begin{eqnarray*}
\frac{F(v)^2}{\dvol}
&=&
\left|\left(W^+ + \frac{s}{12}\right)(v)\right|^2
-
\left| \Ric_0(v)\right|^2,\\
&&\\
&=&
\langle \mathcal D(v), v \rangle.
\end{eqnarray*}
This proves the statements concerning definiteness of $\mathcal D$ and the orientations.

In the first case, to check the sign of the connection, consider the path of maps 
$$
F_t = \left(W^+ +\frac{s}{12}\right) \oplus t \Ric_0
$$
for $t\in [0,1]$. When the inequality holds all the maps $F_t$ have definite image and so carry a sign which, by continuity, must be independent of $t$. For $t=1$ this is the sign we seek, but for $t=0$ it is simply the sign of $\det \left(W^+ + \frac{s}{12}\right)$.

In the second case, to check the sign run a similar argument with the path
$$
F_t = t\left(W^+ + \frac{s}{12}\right)\oplus \Ric_0.
$$
(Here the sign of $\det \Ric_0$ makes sense as $\Ric_0 \colon \Lambda^+ \to \Lambda^-$ is an invertible map between oriented bundles.)
\end{proof} 

Of course, we could have worked throughout with the bundle $\Lambda^-$ giving similar conditions involving $W^-$ in place of $W^+$. 

It follows from Proposition \ref{half Hitchin-Thorpe}, that compact manifolds with metrics for which $\mathcal D>0$  have $2\chi + 3\tau>0$, whilst when $\mathcal D<0$, $2\chi-3\tau>0$ (as the definite connection induces the opposite orientation). In fact, when $\mathcal D<0$ we can say more: $\Ric_0$ gives an isomorphism $\Lambda^+\to\Lambda^-$; since $p_1(\Lambda^+)=2\chi+3\tau$ whilst $p_1(\Lambda^-) = -2\chi + 3\tau$ it follows that $\chi=0$ and $\tau<0$.

\begin{proposition}\label{Riemannian Hitchin-Thorpe}
For a compact Riemannian four-manifold, $\mathcal D>0$ implies $2\chi+3\tau>0$ whilst $\mathcal D<0$ implies $\chi=0$ and $\tau<0$.
\end{proposition}

In Riemannian geometry, these Hitchin--Thorpe type inequalities are usually proved via the Chern--Weil integral formula
\begin{equation}\label{2chi+3tau}
2\chi + 3\tau
=
\frac{1}{2\pi^2}\int_X\left(
|W^+|^2 +\frac{s^2}{48} - |\Ric_0|^2
\right).
\end{equation}
See, for example, \cite{besse}; note that the coefficient of $|\Ric_0|^2$ used here is different because we are using the norm of the trace-free Ricci \emph{operator} $\Ric_0\colon \Lambda^+\to\Lambda^-$, this is one-quarter the size of the norm of the trace-free Ricci \emph{curvature} $\ric_0\colon TX \to TX$. (The two are related by equation (\ref{trace-free Ricci operator}).)

To see the inequality $2\chi+3\tau>0$ from this more traditional perspective, note that $\mathcal D$ is essentially the ``$F^2$'' in the Chern--Weil representative $\tr F^2$ for $p_1(\Lambda^+)$. In particular,
$$
\tr \mathcal D
=
|W^+|^2 + \frac{s^2}{48} - |\Ric_0|^2
$$
is the integrand in (\ref{2chi+3tau}). Of course, $\mathcal D>0$ implies $\tr \mathcal D>0$ and hence $2\chi+3\tau>0$.

The inequality given by requiring $\mathcal D$ to be definite determines an open cone in the space of Riemannian curvature tensors. This cone has six components, in contrast to the case of a general definite connections where there are only four components. This is because the components are separated not only by the induced orientation and sign of the definite connection but also, when $\mathcal D >0$, by the type of the non-degenerate symmetric form $(W^++\frac{s}{12})$; if the determinant is positive we still have the possibility that the signature is $3$ or $-1$, if the determinant is negative we still have the possibility that the signature is $-3$ or $1$. 

At this stage we should stress that the only complete metrics satisfying Proposition \ref{Riemannian curvature inequality} that we have found so far lie in two of these six components: they all have $\mathcal D>0$ and, moreover, in all cases $W^++\frac{s}{12}$ is definite.

\subsection{Anti-self-dual solutions}\label{asd solutions}

\subsubsection{Anti-self-dual Einstein metrics}

Given the appearance of $W^+$ in Proposition \ref{Riemannian curvature inequality} it is natural to look for solutions with $W^+=0$. Before describing general anti-self-dual solutions and the closely related work of Gauduchon \cite{gauduchon} we first discuss, at long last, some examples. When both $W^+$ and $\Ric_0$ vanish then $\mathcal D >0$ is equivalent to $s \neq 0$.

\begin{corollary}
For an anti-self-dual Einstein four-manifold with non-zero scalar curvature, the Levi--Civita connection on $\Lambda^+$ is definite. It induces the same orientation, and carries sign equal to that of the scalar curvature.
\end{corollary}

The round metric on $S^4$ is Einstein and conformally flat---i.e., anti-self-dual in either orientation---hence gives positive definite connections inducing both orientations. Meanwhile, the Fubini--Study metric on $\C\P^2$ is self-dual and Einstein; after changing the orientation, the metric becomes anti-self-dual and so gives a positive definite connection inducing the non-standard orientation. In both cases these metrics lead to K\"ahler Fano manifolds: $S^4$ gives $\C\P^3$ whilst $\C\P^2$ gives the complete flag $F(\C^3)$. (It is well known that the sphere bundles are diffeomorphic to the manifolds claimed; we will justify in \S\ref{anti-self-dual solutions} that they are symplectomorphic.) In fact, a theorem of Hitchin \cite{hitchin2} says that these are the only compact examples of anti-self-dual Einstein metrics with positive scalar curvature and so there are no new compact symplectic Fano manifolds to be found this way. 

On the other hand, the corollary does lead to interesting symplectic Calabi--Yau manifolds. For example, hyperbolic four-manifolds are Einstein and conformally flat so have negative definite connections inducing either orientation.  Meanwhile, complex-hyperbolic manifolds are Einstein and self-dual and so admit negative definite connections inducing the non-standard orientation. The symplectic manifolds arising this way were known to Reznikov \cite{reznikov} and also Davidov--Grantcharov--Mu\v{s}karov \cite{davidov} (although seemingly not that they are Calabi--Yau). 

These symplectic Calabi--Yau manifolds have certain properties that K\"ahler Calabi--Yaus do not.

\begin{lemma}
~
\begin{enumerate}
\item
For a symplectic Calabi--Yau arising from a hyperbolic four-manifold $X$, $c_2 \cdot [\omega] = -2 \chi(X)$. 
\item
For a symplectic Calabi--Yau arising from a complex-hyperbolic four-manifold, $c_2 \cdot [\omega] = 0$. 
\end{enumerate}
\end{lemma}

\begin{proof}
Recall that Corollary \ref{Chern numbers} gives $c_2 \cdot [\omega] = -2(\chi+3\tau)$. For part 1, note that hyperbolic four-manifolds, like all conformally flat four-manifolds, have zero signature. 

For part 2, note that complex-hyperbolic four-manifolds are precisely those K\"ahler surfaces with $c_1^2= 3c_2$; in other words, with respect to the \emph{complex} orientation, $3\tau = \chi$. However, as the metric is self-dual in the complex orientation, the orientation induced on $X$ by the negative definite connection is opposite to the complex one and so, in the formula from Corollary \ref{Chern numbers}, we must use $3\tau = - \chi$.
\end{proof}

To put this result in context, note that for a \emph{K\"ahler} Calabi--Yau threefold, $c_2\cdot[\omega] \geq 0$; moreover, $c_2 \cdot [\omega]=0$ only when the manifold admits a flat metric. This is because, when using the Ricci flat metric, the Chern--Weil formula shows that this quantity is precisely the $L^2$-norm of the curvature tensor. In \cite{smith-thomas-yau}, Smith--Thomas--Yau were led to ask whether this quantity could be negative for a symplectic Calabi--Yau. The symplectic manifolds associated to hyperbolic manifolds show that it can.

As is mentioned in \cite{smith-thomas-yau,gross}, the positivity of $c_2\cdot[\omega]$ in the K\"ahler setting has an interpretation in terms of the Lagrangian torus fibrations which appear in the SYZ conjectural description of mirror symmetry \cite{syz}. If a symplectic manifold $Z$ has a suitably nice Lagrangian torus fibration then it must have both $c_1=0$ and $c_2\cdot[\omega]\geq 0$. For a rigorous discussion (including a description of ``suitably nice'') we refer to \cite{gross}. The idea is that away from the singular locus $\Sigma$, the tangent bundle $TZ$ is symplectically trivial and so the Chern classes of $Z$ are concentrated on $\Sigma$. Under certain hypotheses $\Sigma$ is actually a symplectic surface and so $c_1=0$ whilst $c_2\cdot [\omega] \geq0$. This leads to the natural question: do the symplectic Calabi-Yaus coming from hyperbolic four-manifolds admit any kind of Lagrangian torus fibrations? 

\subsubsection{Fake projective planes}\label{fake projective planes}

One class of complex-hyperbolic manifolds deserves a special mention. A \emph{fake projective plane} is a compact complex surface with the same Betti numbers as $\C\P^2$. Since the first example of such a surface---produced by Mumford \cite{mumford}---much research has focused on determining all such surfaces. The problem remains open, but to date there are several known examples. The article \cite{prasad-yeung} contains the best results known to-date; for more details on examples, see this paper and the references therein.

For our purposes, the important fact is that fake projective planes automatically admit complex-hyperbolic metrics and hence give symplectic Calabi--Yau manifolds with $b_1=0=b_3$. As explained in Lemma \ref{holomorphic volume form} this is impossible for a K\"ahler Calabi--Yau threefold due to the presence of a holomorphic volume form. In \cite{smith-thomas-yau} Smith--Thomas--Yau wondered whether this condition could hold for a symplectic Calabi--Yau and we see here that in fact it can.

\subsubsection{Anti-self-dual metrics with definite Ricci operator}\label{anti-self-dual solutions}

We return now to the discussion of anti-self-dual metrics which satisfy Proposition \ref{Riemannian curvature inequality}. When the metric is anti-self-dual, the twistor space $Z$ is, in a natural way, a complex threefold. This is the starting point of the famous Penrose correspondence, developed in the Riemannian context by Atiyah--Hitchin--Singer \cite{atiyah-hitchin-singer}. Using this complex perspective, Gauduchon \cite{gauduchon} gives an alternative definition of the 2-form $\omega$ on $Z$. 

In the anti-self-dual setting, it turns out that the vertical tangent bundle $V\to Z$ is naturally a holomorphic Hermitian line bundle. Gauduchon considers the Chern connection of this line bundle, which is the same as our connection $\nabla^V$. As it is the curvature of a holomorphic line bundle, it follows that $\omega$ is a $(1,1)$-form and Gauduchon studies the condition that it be a K\"ahler metric. This is not quite the same thing as $\omega$ being symplectic, however. The distinction is illustrated by the form on $\C^3$
\begin{equation}\label{symplectic non-Kahler form}
\frac{i}{2}\left(
\diff z_1 \wedge \diff \bar z_1
-\diff z_2 \wedge \diff \bar z_2
-\diff z_3 \wedge \diff \bar z_3\right).
\end{equation}
This real $(1,1)$-form is symplectic but only becomes K\"ahler after the complex structure is reversed in the $(z_2,z_3)$-plane. 

To describe Gauduchon's result we first need to introduce the Ricci operator; this is the Riemann curvature operator minus its Weyl part. Explicitly, as an endomorphism of $\Lambda^2=\Lambda^+\oplus\Lambda^-$, it is given by
$$
\Ric
=
\left(\begin{array}{cc}
\frac{s}{12} & \Ric_0^*\\
&\\
\Ric_0 & \frac{s}{12}
\end{array}\right).
$$
Gauduchon proved:
\begin{theorem}[Gauduchon, \cite{gauduchon}]
Given an anti-self-dual metric, the form $\omega$ is K\"ahler if and only if the Ricci operator is positive definite
\end{theorem}

A theorem of Hitchin \cite{hitchin2} guarantees that the only compact K\"ahler twistor spaces are those belonging to the standard conformal structures on $S^4$ and $\overline{\C\P}^2$. It follows that a compact anti-self-dual metric with positive definite Ricci operator is conformally equivalent to one of these two spaces. Moreover, the symplectic Fanos arising from such metrics are precisely those K\"ahler twistor spaces, namely $\C\P^3$ and the complete flag $F(\C^3)$.  

In this article we are interested in the more general question of whether or not $\omega$ is symplectic. To this end we obtain:

\begin{proposition}
An anti-self-dual metric has $\mathcal D>0$ if and only if it has definite Ricci operator. When this happens, the sign of the corresponding definite connection agrees with the sign of the Ricci operator.
\end{proposition}
\begin{proof}
When $W^+=0$, positivity of $\mathcal D$ is equivalent to the positivity of the operator
$$
\frac{s^2}{144} - \Ric_0^*\Ric_0 \colon \Lambda^+ \to \Lambda^+
$$
which is in turn equivalent to the definiteness of the Ricci operator. When this happens, the sign of the definite connection agrees with the sign of the scalar curvature which agrees with the sign of the Ricci operator.
\end{proof}

As we will see in \S\ref{twistor spaces}, the symplectic structures corresponding to anti-self-dual metrics with negative definite Ricci operator fail to be K\"ahler for precisely the same reason as the form (\ref{symplectic non-Kahler form}) above.

Gauduchon's result shows that there are no interesting Fano manifolds to be found directly via anti-self-dual metrics. On the other hand, LeBrun--Nayatani--Nitta \cite{lebrun-nayatani-nitta}, answering a question posed by Gauduchon, have shown that $\overline{\C\P}^2\sharp\overline{\C\P}^2$ admits an anti-self-dual metric with \emph{non-negative} Ricci operator; moreover, it is positive on a dense set. This metric lies on the boundary of the inequality $\mathcal D>0$; it would be interesting to know whether or not it can be perturbed to give a positive definite connection. 

We close this discussion of the Ricci operator with a warning to readers who have not met it before: definiteness of the Ricci operator is stronger than definiteness of the Ricci curvature thought of as an endomorphism $TX\to TX$. Indeed, if $\lambda_1, \lambda_2, \lambda_3, \lambda_4$ denote the eigenvalues of the Ricci curvature when considered as an endomorphism of $TX$; then the eigenvalues of the Ricci operator are readily seen to be $\frac{1}{2}(\lambda_i+\lambda_j)-\frac{s}{6}$ for $i<j$. The Ricci operator is positive definite when 
$\lambda_i+\lambda_j>\frac{s}{3}$ from which it follows that in fact $\frac{2s}{3}>\lambda_i+\lambda_j>\frac{s}{3}$. These inequalities imply that all $\lambda_i>0$, but the additional pinching they describe is even more restrictive and so positive Ricci operator is far stronger than positive Ricci curvature in the traditional sense. Likewise, negative Ricci operator is far stronger than negative Ricci curvature.

\subsection{Metrics with pointwise-pinched sectional curvatures}

Recall that a metric on $X$ is said to have pointwise $K$-pinched sectional curvatures if all sectional curvatures have the same sign and, at any point $p\in X$, 
$$
\frac{\min_p \sec}{\max_p \sec} > K,
$$
(where $\min_p \sec$ denotes the minimum sectional curvature at $p$ and $\max_p\sec$ the maximum at $p$). Now, the curvature inequality $\mathcal D > 0$ is open and is satisfied by the four-sphere which has all sectional curvatures equal to one. It follows that there must be a constant $K<1$ such that any metric with positive sectional curvatures pointwise $K$-pinched also has $\mathcal D>0$. Similarly, by considering metrics close to the hyperbolic metric, there must be a constant $K<1$ such that any metric with negative sectional curvatures pointwise $K$-pinched also has $\mathcal D>0$. (A similar observation was made by Reznikov \cite{reznikov}.) 

This main result of this section is that this is true for $K=2/5$ and that this choice of $K$ is optimal (in both positive or negative cases).

\begin{theorem}\label{pinching theorem}~
\begin{enumerate}
\item
Any metric with positive sectional curvatures pointwise $2/5$-pinched induces positive definite connections on both $\Lambda^+$ and $\Lambda^-$.
\item
Any metric with negative sectional curvatures pointwise $2/5$-pinched induces negative definite connections on both $\Lambda^+$ and $\Lambda^-$.
\end{enumerate}
Moreover, in each case 2/5 is the optimal pinching constant for which this result remains true.
\end{theorem}

\begin{proof}
We consider only the positive case since the proof in the negative case is identical. Let $\mathcal R$ denote the set of all algebraic curvature tensors, i.e., all tensors in $S^2(\Lambda^2 \R^4)$ satisfying the algebraic Bianchi identity. Given $K$, let $P_K \subset \mathcal R$ denote those curvature tensors with positive $K$-pinched sectional curvatures. Since sectional curvatures are insensitive to orientation, if all tensors in $P_K$ have $\mathcal D>0$, they also solve the same inequality with the opposite orientation (in other words, with $W^+$ replaced by $W^-$ in definition (\ref{main operator}) and so inducing definite connections on both $\Lambda^+$ and $\Lambda^-$). Let $B\subset \mathcal R$ denote all curvature tensors which have $\mathcal D>0$ in \emph{both} orientations. Each $P_K$ is connected, its boundary varies continuously with $K$ and, for $K_1<K_2$, $P_{K_2}\subset P_{K_1}$, so to prove our claim we have to show that as we decrease $K$ from 1, the first time the closure of $P_K$ meets the boundary of $B$ is when $K= 2/5$. Note that $B$ has several components according to the types of the non-degenerate operators $W^\pm +s/12$. As we are interested in the component containing the four-sphere, we focus on the component of $B$ in which $W^\pm +s/12$ are both positive definite.

We write a curvature tensor $R \in\mathcal R$ with respect to the decomposition  $\Lambda^2=\Lambda^+ \oplus \Lambda^-$ as
$$
R = \left(
\begin{array}{cc}
A & B^*\\
B & C
\end{array}
\right).
$$
Sectional curvatures are found by evaluating $R$ on decomposable elements of $\Lambda^2$ which are precisely those of the form $u+v$ where $u\in \Lambda^+$ and $v\in \Lambda^-$ have the same length. When $|u| = |v| = 1$, the corresponding sectional curvature is
$$
\sec(u+v) = (Au,u) + 2(Bu,v) + (Cv,v).
$$
 
Assume now that $R\in\del B$ with both $A$ and $C$ positive semi-definite; we want to show that $\frac{\min\sec(R)}{\max \sec (R)}\leq \frac{2}{5}$. By scaling we can assume that $s =12$ so that $\tr A = 3  = \tr C$. Denote the highest and lowest eigenvalues of $A$ and $C$ by $1+a_1$, $1 - a_2$ and $1+c_1$, $1-c_2$ respectively. Since $A$ and $C$ are positive semi-definite we have
\begin{equation}\label{ineq1}
2 \geq 2a_2 \geq a_1\geq \frac{a_2}{2} \geq 0,\quad
2 \geq 2c_2 \geq c_1\geq \frac{c_2}{2} \geq 0
\end{equation}

Next we claim the following inequalities for $R$:
\begin{equation}\label{ineq2}
\max \sec(R) \geq 2 +a_1 + c_1,\quad \min \sec(R)\leq 2-a_2-c_2.
\end{equation}
\begin{equation}\label{ineq3}
\min \sec(R) \leq c_1 + a_2\quad \mathrm{or} \quad
\min \sec(R) \leq a_1 + c_2.
\end{equation}
The inequalities (\ref{ineq2}) hold for any $R$ with $s=12$. To prove the first of them, let $u_1$ and $v_1$ be unit-length eigenvectors of $A$ and $C$ with eigenvalue $1+a_1$ and $1+c_1$ respectively. Then one of $\sec(u_1 \pm v_1)$ is greater than  $2+a_1+c_1$. The second inequality in (\ref{ineq2}) is proved similarly. 

To prove (\ref{ineq3}) note that, since $R$ is on the boundary of $B$, there is either a unit vector $u$ such that $|Au| = |Bu|$ or a unit vector $v$ such that $|Cv|=|B^*v|$. We assume the first occurs, the second case leads in an identical way to the second inequality in (\ref{ineq3}). Let $v = -\frac{Bu}{|Bu|}$. Note that $|Au|\geq 1 - a_2$ whilst $|Cv| \leq 1+c_1$, so
\begin{eqnarray*}
\sec(u+v)& \leq		& |Au| - 2|Bu| + 1 + c_1,\\
			& 		=	& -|Au| + 1 +c_1,\\
			& \leq 	& c_1 + a_2.
\end{eqnarray*}

With inequalities (\ref{ineq1}) (\ref{ineq2}) and (\ref{ineq3}) in hand we proceed by contradiction and suppose that  
$
\frac{\max \sec(R)}{\min \sec(R)} <5/2
$.
From (\ref{ineq2}) and (\ref{ineq3}) (and assuming the first case in (\ref{ineq3}) holds, the second case being similar) we find
\begin{eqnarray*}
\frac{5}{2}(2-a_2-c_2) &>& 2 + a_1 + c+1,\\
\frac{5}{2}(c_1+a_2) &> &2 + a_1 + c_1.
\end{eqnarray*}
Using $c_2\geq c_1/2$ we deduce
\begin{eqnarray*}
3 - \frac{5}{2}a_2 &>& a_1 + \frac{9}{4}c_1,\\
\frac{3}{2}c_1 + \frac{5}{2}a_2 &>& 2 +a_1.
\end{eqnarray*}
Rewriting this as
$$
3 - \frac{5}{2}a_2 - a_1 
> 
\frac{9}{4}c_1
>
3+\frac{3}{2}a_1-\frac{15}{4}a_2,
$$
hence 
$$
\frac{5}{4}a_2 - \frac{5}{2}a_1 >0
$$
which contradicts $2a_1\geq a_2$.

This establishes that $P_{2/5}$ is contained in $B$. To prove that $2/5$ is the optimal pinching constant, consider the curvature tensor with $B=0$, $C=1$ and $A=\diag(3/2, 1, 0)$. This lies in the closure of $P_{2/5}$ and the boundary of $B$. 
\end{proof}

\begin{remark}\label{pinched metrics tame}
In fact, metrics with sectional curvatures pointwise $2/5$-pinched satisfy stronger inequalities than just $\mathcal D>0$. For example, in the case of positive sectional curvatures, for unit length $u\in \Lambda^+$, $v\in \Lambda^-$, both the following hold:
\begin{eqnarray*}
(Au, u) &>& |Bu|,\\
(Cv, v) &>& |B^*v|.
\end{eqnarray*}
(These stronger inequalities will be important later on when we relate the symplectic form to twistor geometry---see \S\ref{twistor spaces}.)

The proof is identical to that above. Letting $B'$ denote those curvature tenors solving the above inequalities, we simply have to check that for $R \in \del B'$ with both $A$ and $C$ positive semi-definite, one of the inequalities (\ref{ineq3}) still holds. This time we assume that there is a unit vector $u$ such that $(Au,u) = |Bu|$ and consider the sectional curvature corresponding to $u+v$ where, as before, $v = -Bu/|Bu|$; the assumption is still enough to deduce that $\sec(u+v) < c_1 + a_2$ and the rest of the proof goes through unchanged. Similarly, metrics with negative sectional curvatures pointwise $2/5$-pinched satisfy $|(Au,u)|>|Bu|$ and $|(Cv,v)| > |B^*v|$ for all unit length $u\in \Lambda^+$, $v\in \Lambda^-$.
\end{remark}

As with anti-self-dual solutions, it is not possible to find any new symplectic Fanos via pinched metrics: a  theorem of Chen \cite{chen} says that the only orientable compact four-manifold which admits a metric with positive sectional curvatures pointwise $1/4$-pinched is $S^4$.

We close this section by recalling a result due to Ville \cite{ville} that a compact oriented Riemannian four-manifold with negative sectional curvatures strictly $1/4$-pinched must have $3|\tau|<\chi$. Note that for Ville's result, the pinching must be \emph{global}, and not just pointwise, i.e., all sectional curvatures satisfy $-1<\sec < -1/4$. It follows that those metrics with negative sectional curvatures globally $2/5$-pinched give more examples of symplectic Calabi-Yaus with $c_2\cdot[\omega] = -2(\chi+3\tau)<0$.

\subsection{Gromov--Thurston metrics}
\label{GT metrics}

In \cite{gromov-thurston} Gromov and Thurston construct interesting examples of negatively curved metrics via ramified coverings of hyperbolic manifolds. In one direction, they show that for any small $\epsilon>0$ there is a Riemannian four-manifold with negative sectional curvatures pointwise $(1-\epsilon$)-pinched but which admits no hyperbolic metric. It follows from the previous section that these manifolds carry negative definite connections on both $\Lambda^+$ and $\Lambda^-$ providing examples of such connections on manifolds which are not quotients of symmetric spaces.

In another direction they use the same construction to produce, for any $K$, Riemannian four-manifolds with negative sectional curvatures but which admit no pointwise $K$-pinched metric. In this section we will show that these metrics on ``unpinchable'' manifolds also have $\mathcal D >0$ in both orientations.

We begin by reviewing the Gromov--Thurston construction. Take hyperbolic space $H^n$ and a totally geodesic subspace $H^{n-2}\subset H^n$. The hyperbolic metric on $H^n$ can be written as
$$
\diff r^2 + \sinh^2(r)\diff \theta^2 + \cosh^2(r) g_{H^{n-2}}
$$
where $(r,\theta)$ are polar coordinates in the normal directions to $H^{n-2}$ and $g_{H^{n-2}}$ is the hyperbolic metric on $H^{n-2}$. Let $\widetilde H^n$ denote the $k$-fold cover of $H^n$ ramified along $H^{n-2}$. Of course, the hyperbolic metric becomes singular when lifted to $\widetilde H^n$ but this singularity can be smoothed out by replacing the coefficient $\sinh^2(r)$ of $\diff \theta^2$ with a function $\sigma(r)^2$ where $\sigma(r) = k \sinh(r)$ for small $r$. A judicious choice of $\sigma$ ensures that the sectional curvatures of $\widetilde H^n$ remain negative and that their pointwise pinching can be controlled. Moreover, taking $\sigma(r) = \sinh(r)$ for $r$ greater than some $r_0$ means that outside of some neighbourhood of the ramification locus the metric on $\widetilde H^n$ is hyperbolic. 

To produce compact examples, Gromov--Thurston construct closed hyperbolic manifolds each containing a totally geodesic codimension 2 submanifold with normal injectivity radius less than $r_0$. The metric on $\widetilde H^n$ then descends to give metrics on the covers of these hyperbolic manifolds ramified along the totally geodesic submanifolds.

For our purposes all that matters is the curvature tensor of $\widetilde H^n$ and Gromov--Thurston compute this explicitly. Let $x_1, \ldots x_{n-2}$ denote coordinates on $H^{n-2}$, then the curvature tensor is diagonal for the obvious basis of $\Lambda^2$ coming from wedging pairs from the basis $\del_{x_i}$, $\del_r$, $\del_\theta$ and the diagonal entries are:
\begin{eqnarray*}
\Rm(\del_r, \del_{x_i}, \del_r, \del_{x_i})
&=&
-1, \\
\Rm(\del_{x_i}, \del_{x_j}, \del_{x_i}, \del_{x_j})
&=&
-1 
\quad\mathrm{for}\  i\neq j,\\
\Rm(\del_{x_i}, \del_\theta, \del_{x_i}, \del_\theta)
&=&
-\frac{\sigma'(r)}{\sigma(r)}\tanh(r),\\
\Rm(\del_r,\del_\theta, \del_r,\del_\theta)
&=&
-\frac{\sigma''(r)}{\sigma(r)}.
\end{eqnarray*}
Note that there is a misprint in the published paper which has the sign wrong in the third of these sectional curvatures. The point is that $\sigma$ is chosen so that $\sigma', \sigma'' >0$, implying  that $\Rm$ is negative definite.

We are interested in the case $n=4$. The above description of $\Rm$ is given using the standard basis of $\Lambda^2$ via decomposable forms built from the orthogonal coframe $\diff x_1, \diff x_2, \diff r, \diff \theta$. If we switch to the standard basis of $\Lambda^2$ using self-dual and anti-self-dual forms built from this coframe we get a curvature operator
$$
\Rm=
\left(
\begin{array}{cc}
A & B^*\\
B & C
\end{array}
\right)
$$
in which $A$, $B$ and $C$ are all diagonal. It is now straightforward to check that negative definiteness of $\Rm$ implies that $|A(u)|>|B(u)|$ and $|C(v)|>|B^*(v)|$ for all $u\in\Lambda^+$ and $v\in\Lambda^-$. In other words the Gromov--Thurston metrics satisfy $\mathcal D>0$ in both orientations. 

\subsection{Consequences of the Riemannian inequality}

As we have seen, any compact four-manifold admitting a definite connection has $2\chi + 3\tau >0$. In the case of definite connections coming from Riemannian metrics---and with an additional curvature assumption---one can conclude far more, however. The condition $\mathcal D>0$ implies that $W^+ + s/12$ is an isomorphism $\Lambda^+ \to \Lambda^+$; in this section we will assume that, in addition, it is positive definite. 

The first consequence comes from a standard Bochner argument giving $b_+ =0$. Indeed, letting $A= W^++s/12$, the following is well-known. 

\begin{lemma}
Suppose that the sum of the lowest two eigenvalues of $A$ is positive. Then $b_+ = 0$.
\end{lemma}

\begin{proof}
The Weitzenbock formula for self-dual 2-forms asserts
\begin{eqnarray*}
\Delta
&=&
\nabla^* \nabla
- 2W^+
+
\frac{s}{3},
\\
&=&
\nabla^*\nabla
+
2(\tr A - A).
\end{eqnarray*}
Applying this to a harmonic self-dual 2-form $\alpha$ and integrating against $\alpha$ gives
$$
0
=
\|\nabla \alpha \|^2 
+ 
2 \int_M\langle \tr (A - A)\alpha, \alpha\rangle.
$$
By assumption, $\tr A - A$ is a positive operator. Hence both terms on the right hand side vanish and $\alpha=0$
\end{proof}

Since $\mathcal D>0$ implies $2\chi+3\tau>0$ and $A>0$ implies $b_+=0$, together they give $b_-=b_2< 4(1-b_1)$. Hence $b_1=0$ and $b_-=b_2 < 4$. As Proposition \ref{ricci positive} below shows, the conditions $\mathcal D>0$ and $A>0$ together also imply that the Ricci curvature is positive, which in turn gives even more control over the fundamental group. We begin the proof of this fact by recalling a lemma from linear algebra. 
% The proof of this is in a comment at the end of the document

\begin{lemma}
Let $A, B \colon \R^n \to \R^n$ be two self-adjoint endomorphisms with $|A(v)|>|B(v)|$ for all $v\neq 0$. Let $a_1, \ldots , a_n$ denote the eigenvalues of $A$ and $b_1, \ldots, b_n$ the eigenvalues of $B$. Then $\sum |a_i| > \sum |b_i|$.
\end{lemma}

\begin{proposition}\label{ricci positive}
A oriented Riemannian four-manifold with $\mathcal D>0$ and $W^++s/12>0$ has positive Ricci curvature.
\end{proposition}

\begin{proof}
Let $\lambda_1\geq \lambda_2\geq \lambda_3\geq \lambda_4$ be the eigenvalues of the Ricci curvature $\ric\colon T^*X\to T^*X$. Let $\alpha_1,\alpha_2,\alpha_3, \alpha_4$ denote an oriented orthonormal basis of eigenvectors. The following triples
$$
\{ \alpha_1 \wedge \alpha_2 \pm \alpha_3\wedge \alpha_4,\
\alpha_1\wedge\alpha_3 \mp \alpha_2 \wedge \alpha_4,\
\alpha_1\wedge \alpha_4 \pm \alpha_2\wedge \alpha_4 \}
$$
give orthogonal bases in $\Lambda^\pm$ respectively. Using these bases to identify $\Lambda^\pm\cong \R^3$ we will apply the linear algebra lemma to the operators $A=W^++s/12$ and $B=\Ric_0$. 

By assumption, all eigenvalues of $A$ are positive and so $\sum|a_i|= \tr A = s/4$. Moreover, in these bases the trace free Ricci operator $\Ric_0 \colon \Lambda^+ \to \Lambda^-$ is diagonal (in particular self-adjoint) and has eigenvalues
$$
\frac{1}{4}(\lambda_1 + \lambda_2 - \lambda_3 - \lambda_4)
\geq
\frac{1}{4}(\lambda_1 - \lambda_2 + \lambda_3 - \lambda_4)
\geq
\frac{1}{4}(\lambda_1 - \lambda_2 - \lambda_3 + \lambda_4).
$$

Since the $\lambda_i$ are written in descending order, the first two eigenvalues here are non-negative. We split in to two cases according to the sign of the third eigenvalue of $\Ric_0$.

\begin{enumerate}
\item
We assume that $\lambda_1 - \lambda_2 - \lambda_3 + \lambda_4 \geq0$. The inequality $s/4 > \sum |b_i|$ reads $s > 3\lambda_1 - \lambda_2 -\lambda_3 - \lambda_4$ which in turn implies $\lambda_2 + \lambda_3 > \lambda_1 - \lambda_4$. Combining this with the assumption at the start of this case we see that $\lambda_4>0$ and so the Ricci curvature is positive.

\item
We assume that $\lambda_1 - \lambda_2 - \lambda_3 + \lambda_4 <0$. The inequality $s/4 > \sum |b_i|$ reads $s > \lambda_1 + \lambda_2 + \lambda_3 - 3\lambda_4$ which it turn implies directly that $\lambda_4>0$. Hence, again, the Ricci curvature is positive.
\end{enumerate}
\end{proof}

Combining this with the information from the Bochner argument shows that the conditions $\mathcal D >0$ and $A>0$ strongly restrict the homeomorphism type of the underlying manifold.

\begin{theorem}
Let $X$ be a compact oriented Riemannian four-manifold with $\mathcal D>0$ and $A>0$. Then $X$ is homeomorphic to $n\overline {\C\P}^2$ for some $0\leq n\leq 3$.
\end{theorem}

\begin{proof}
By Proposition \ref{ricci positive}, $X$ has positive Ricci curvature and therefore finite fundamental group. Hence its universal cover $\widetilde X$ is compact, and also has $\mathcal D> 0$ and $A>0$. From the above discussion we see that $\widetilde X$ has $b_+=0$ and $b_2 \leq 3$. Now, thanks to the celebrated theorems of Freedman \cite{freedman} and Donaldson \cite{donaldson}, we conclude that $\widetilde X$ is homeomorphic to $n\overline {\C\P}^2$ for some $0\leq n \leq 3$.

There are various ways to rule out oriented quotients of these manifolds by finite groups. For example, the Lefschetz fixed point theorem says that any orientation preserving diffeomorphism of $n\overline{\C\P}^2$ with no fixed points induces an isomorphism on $H^2$ with trace $-2$, hence none exist for $n=0,1$. To rule out quotients for $n=2,3$ one can use the fact that Euler characteristic and signature are multiplicative in finite unramified covers. For $3\overline{\C\P}^2$, $2\chi + 3\tau = 1$; as this is not divisible by any integer, there is no fixed-point-free orientation preserving diffeomorphism. For $2\overline{\C\P}^2$, $2\chi + 3\tau=2$ and so any such diffeomorphism would have order 2. But then the quotient would have $\chi =2$, as well as $\pi_1 = \Z_2$ and so $H^1 = 0$. Hence the quotient would have $H^2=0$ which  contradicts with the fact that it would also have $\tau = -1$.
\end{proof}

\section{Definite connections and minimal surfaces}\label{twistor spaces}

This section explains how, in certain situations, the existence of a symplectic form arising from a Riemannian four-manifold as above leads to results about minimal surfaces. To set the scene we begin by reviewing an analogous situation in dimension three. 

\subsection{Two classical results in dimension three}

Let $M$ be a compact negatively curved three-manifold. There are two classical results concerning minimal immersions $\Sigma \to M$. Firstly, $\chi(\Sigma)<0$ and secondly, the set of all minimal immersions of $\Sigma$ is weakly compact, i.e., any sequence of such immersions has a subsequence which converges to a minimal immersion possibly with branch points. (We review the definition of branch points later when we consider the four-dimensional case in detail.)

The standard proof of these results begins with Gauss' equation. For a minimal immersion $\Sigma \to M$ with second fundamental form $B$,
$$
K_\Sigma + |B|^2 = \sec(T\Sigma),
$$
where $K_\Sigma$ is the Gauss curvature of the induced metric on $\Sigma$, whilst $\sec(T\Sigma)$ is the sectional curvature of the ambient metric on the plane $T\Sigma \subset TM$. When $M$ is negatively curved there are some immediate consequences: $K_\Sigma$ is everywhere negative and hence $\chi(\Sigma)<0$; moreover, $K_\Sigma < \max \sec (M) <0$ which gives an a priori area bound on the image of $\Sigma$; finally, integrating Gauss' equation gives an a priori bound on the $L^2$-norm of the second fundamental form. These two a priori bounds are the key to an analytic proof of the compactness result mentioned above. One way to proceed is to use a covering argument in the style of Sacks--Uhlenbeck \cite{sacks-uhlenbeck} in a manner which is now standard. (See the article of Choi--Schoen \cite{choi-schoen} for an example of this in the context of minimal surfaces.)

There is an alternative way to proceed, however, using symplectic geometry and it is this approach which we will extend to dimension four. The key observation, due originally to Weierstrass, is the relationship between minimal surfaces and holomorphic curves. We can describe this as follows.

To begin with we need no curvature assumptions on $M$. Let $\pi\colon Y\to M$ denote the unit tangent bundle of $M$. We recall the definition of the contact distribution $W$ in $Y$. The Levi--Civita connection gives a splitting
$$
T_pY \cong V_p\oplus T_{\pi(p)M}
$$
where $V_p$ is the vertical tangent space at $p$. Let $U_p = \langle p \rangle ^\perp \subset T_{\pi(p)}M$ and define the hyperplane distribution $W$ by 
$$
W_p = V_p \oplus U_p.
$$
The point is the following: assume that $M$ and $\Sigma$ are oriented; then any immersion $f \colon \Sigma \to M$ lifts to a map $\tilde f \colon \Sigma \to Y$ where $\tilde f(\sigma)$ is the unit normal to $f_*(T_\sigma \Sigma)$; by construction, the lift is everywhere tangent to $W$.

The next step is to define an almost complex structure $J$ on $W$. To do this take $J$ on $V_p$ to be the standard complex structure on $S^2$; meanwhile the real 2-plane $U_p$ is oriented by the choice of $p$ as its positive normal, so we can take $J$ to be rotation by $-\pi/2$ on $U_p$. Weierstrass' key observation is that an immersion $f$ is minimal if and only if its lift $\tilde f$ is $J$-holomorphic.

\begin{theorem}[The Weierstrass correspondence]
There is a one-to-one correspondence between branched minimal immersions in $M$ and non-vertical $J$-holomorphic curves in $Y$ (which are necessarily tangent to $W$).
\end{theorem}
(Here, non-vertical means we ignore $J$-holomorphic curves contained in a fibre of $Y\to M$.) Minimal immersions in $M$ correspond to $J$-holomorphic curves in $Y$ which are everywhere tangent to the fibres of $Y$, whilst curves with vertical tangencies give branched immersions downstairs. (Again, the precise definite of a branched immersion is given in the more detailed four-dimensional discussion.)

Now, the Levi--Civita connection on $M$ is an $\SO(3)$-connection and so determines a closed 2-form $\omega$ on $Y$. Of course, $\omega$ can't possibly be symplectic, but it can happen that on $W$ it tames $J$ (i.e., $\omega(Ju,u)>0$ for all non-zero $u\in W$). Unwinding the definitions and using Proposition \ref{decomposition of omega} gives the following result.

\begin{theorem}
The natural 2-form $\omega$ on $Y$ tames $J$ on $W$ if and only if $M$ is negatively curved
\end{theorem}

So, when $M$ is negatively curved, we can apply Gromov's theory of $J$-holomorphic curves \cite{gromov}. (This is usually phrased for symplectic manifolds, but there is no difficulty in adapting the arguments to curves tangent to a symplectic hyperplane distribution as is considered here.) One fundamental result is a $J$-holomorphic curve has positive symplectic area. In our situation, for a minimal immersion $f\colon\Sigma\to M$, the symplectic area of the lift is simply $ -\chi(\Sigma)$ and so this translates into the classical fact that $\chi(\Sigma)<0$. Another fundamental result is Gromov's compactness theorem for sequences of $J$-holomorphic curves with fixed symplectic area. (This result is described below during the discussion of the four-dimensional situation.) In the case of minimal immersions in negatively curved manifolds, this gives the classical compactness result mentioned at the start of the section. (One must work a little to show no bubbling can occur; for example, non-vertical bubbles in $Y$ are ruled out because they would give minimal spheres in $M$ contradicting $\chi<0$.)

In fact, this approach is not so different from the classical one outlined above. Both rely on a topological formula ($\int K_\Sigma \diff A = \chi(\Sigma)$ in the classical case, $-\int \tilde f^*\omega = \chi$ in the symplectic case) and, moreover, one can prove Gromov compactness precisely via a Sachs--Uhlenbeck style covering argument \cite{wolfson}. The statement that $\omega$ tames $J$ is a way to neatly encapsulate the a priori bounds for minimal immersions arising from Gauss' equation. Moreover, it shows how these standard results in minimal surface theory fit into a picture familiar to symplectic geometers. 

We will now extend this symplectic approach to minimal surfaces in four-manifolds. One interesting outcome is that the curvature inequality we find is not merely ``negative curvature''. It is possible to find curvature tensors which satisfy the four-dimensional inequality which have some sectional curvatures positive.

\subsection{Taming twistorial almost complex structures}

The first step is to find the analogue in dimension four of the almost complex structure on $Y$. In fact, there are two possibilities here, leading to two very different classes of surface.

Given an oriented Riemannian four-manifold $X$, the unit sphere bundle $Z$ in $\Lambda^+$---the twistor space of $X$---carries two natural almost complex structures, $J_+$ and $J_-$. Atiyah--Hitchin--Singer \cite{atiyah-hitchin-singer} introduced $J_+$ and showed that it is integrable when $W^+=0$; Eells--Salamon \cite{eells-salamon} introduced $J_-$ which, by contrast with $J_+$, is never integrable. To describe $J_\pm$, begin with the splitting $TZ=V\oplus H$ into vertical and horizontal parts induced by the Levi--Civita connection. Both almost complex structures respect this splitting. A point $p\in Z$ corresponds to a self-dual 2-form and so, via the metric, to a skew endomorphism of $T_{\pi(p)}X$. Since the 2-form is self-dual, when suitably scaled the corresponding endomorphism has square $-1$ and so gives an almost complex structure on $T_{\pi(p)}X \cong H_p$; this defines both $J_\pm$ on the horizontal distribution. On the vertical distribution use the standard complex structure on the sphere to define $J_+$ and the complex structure for the opposite orientation on the sphere to define $J_-$. 
 
As is described in the previous sections, the Riemannian metric also determines a natural closed 2-form $\omega$ on $Z$ and we ask when does $\omega$ tame either of $J_\pm$? (Recall that a 2-form $\omega$ tames an almost complex structure $J$ when it is positive on all complex lines, i.e., $\omega(Ju,u)>0$ for all non-zero $u$.) To answer this question we begin with the following preliminary lemma.

\begin{lemma}
Denote the standard linear complex and symplectic structures on $\R^4$ by $J_0$ and $\omega_0$ respectively and let $\alpha \in \Lambda^-(\R^4)$. Then $\omega_0+\alpha$ tames $J_0$ if and only if $|\alpha|^2<2$.
\end{lemma}

\begin{proof}
Suppose $|\alpha|^2<2$. Now, for any $u\neq0$,
$$
(\omega_0+\alpha)(J_0u, u) = |u|^2 + \alpha(J_0u, u),
$$
whilst $|\alpha(u,v)|^2<\frac{1}{2}|\alpha|^2|u|^2|v|^2$. So $|\alpha(J_0u,u)|<|u|^2$ and hence $(\omega_0+\alpha)(J_0u,u)>0$.

For the converse, suppose $\omega_0+\alpha$ tames $J_0$. Then $(\omega_0+\alpha)^2>0$. But
$$
(\omega_0+\alpha)^2 = (2-|\alpha|^2)\frac{\omega_0^2}{2}.
$$
Hence $|\alpha|^2<2$.
\end{proof}

\begin{theorem}\label{tame}
Let $X$ be an oriented Riemannian four-manifold whose curvature satisfies
\begin{equation}\label{tame inequality}
\left| \left\langle 
\left(W^++\frac{s}{12}\right)v, v
\right\rangle \right|
>
\left| \Ric_0(v)\right|
\end{equation}
for every unit length $v\in \Lambda^+$. Then the induced symplectic form on $Z$ tames $J_+$ if $\det(W^++s/12)>0$, whilst it tames $J_-$ if $\det(W^++s/12) <0$.
\end{theorem}

\begin{proof}
First note that the inequality implies $\mathcal D >0$, so that, as expected, $\omega$ is necessarily symplectic. The vertical part of $\omega$ tames the vertical part of $J_\pm$ with sign according to the sign of the connection, i.e., the sign of $\det (W^+ + s/12)$.

At the point $v\in Z$ the horizontal part of the symplectic form is
$$
(W^++s/12)(v) + \Ric_0(v)
$$
and we are interested in when this tames the almost complex structure $J_v$ determined by $v$. The first term in this expression is a self-dual two-form. Write $\Lambda^+ = \langle v\rangle \oplus \langle v\rangle^\perp$; the 2-forms in $\langle v\rangle^\perp$ vanish on $J_v$-complex lines, so the value of the above form on $J_v$-complex lines is unchanged by projecting the first term onto $v$. In other words, the symplectic form tames if and only if the form 
$$
\left\langle 
\left(W^++\frac{s}{12}\right)v, v
\right\rangle
+
\Ric_0(v)
$$ 
tames $J_v$. The result now follows from the preceding Lemma.
\end{proof}

\begin{remark}
Notice that when inequality (\ref{tame inequality}) is satisfied, $\mathcal D>0$ and $W^++s/12$ is necessarily definite. In particular, the space of curvature tensors satisfying this more restrictive inequality has only two components, compared with the six components of the space of curvature tensors for which $\mathcal D$ is definite. 
\end{remark}

With the exception of the Gromov--Thurston metrics described in \S\ref{GT metrics}, it is relatively straightforward to check that all the examples of metrics with $\mathcal D>0$ that we have given so far in fact satisfy the tighter inequality (\ref{tame inequality}) of Theorem \ref{tame}. For example, for metrics with sectional curvatures pointwise 2/5-pinched this follows from Remark \ref{pinched metrics tame}. For the Gromov--Thurston metrics it is not clear, to us at least, if inequality (\ref{tame inequality}) is satisfied in general.

\subsection{Pseudoholomorphic curves and minimal surfaces}

Having found Riemannian manifolds with tamed twistor spaces, the next step is to consider the $J_\pm$-holomorphic curves. In \cite{eells-salamon} Eells--Salamon make an in-depth study of the $J_\pm$-holomorphic curves and we review some of their findings briefly here.

Given an immersed oriented surface $f\colon\Sigma \to X$, there is a natural lift $\tilde f \colon \Sigma \to Z$ to the twistor space: given an oriented basis $u, v$ for $T_\sigma \Sigma$, take $\tilde f(\sigma)$ to be the unit-length vector in the direction of the self-dual part of $f_*(u)\wedge f_*(v)$. Alternatively, one can think of the fibres of $Z$ as almost complex structures on $X$ compatible with the orientation and metric. There is a unique such almost complex structure on $T_{f(\sigma)}X$ which makes $f_*(T_\sigma \Sigma)$ a complex line with correct orientation and this gives another way to define $\tilde f(\sigma)$. 

Whether or not $\tilde f(\Sigma)$ is a $J_\pm$-holomorphic curve comes down to the second fundamental form $B$ of $f$. The tensor $B$ is a section of $S^2(T^*\Sigma) \otimes N$, where $N$ is the normal bundle. The lift $\tilde f$ determines an almost complex structure on $TX|_{f(\Sigma)}$ which respects the splitting $TX|_{f(\Sigma)} = T\Sigma \oplus N$ so we can think of both $T\Sigma$ and $N$ as complex vector bundles. Using this we can split the space $S^2(T^*\Sigma)\otimes N$ of \emph{real} bilinear symmetric forms with values in $N$ into parts:
$$
\left(S^{2,0}(T^*\Sigma) \otimes_\C N\right)
\oplus
\left(S^{1,1}(T^*\Sigma) \otimes_\C N\right)
\oplus 
\left(S^{0,2}(T^*\Sigma) \otimes_\C N\right)
$$
where the first summand is complex bilinear, the second is the symmetrisation of complex linear in the first argument and anti-linear in the second whilst the third summand is anti-bilinear.

Under this decomposition the second fundamental form splits into parts $B = B^{2,0}+B^{1,1}+B^{0,2}$. Explicitly, if we write $J$ for the almost complex structure on $TX|_{f(\Sigma)}$ corresponding to $\tilde f$,
\begin{eqnarray*}
B^{2,0}(u,v) &=&
\frac{1}{4}\left[
B(u,v) - J B(Ju,v) - J B(u,Jv) - B(Ju,Jv)
\right],\\
&&\\
B^{1,1}(u,v) &=&
\frac{1}{2}\left[
B(u,v) + B(Ju,Jv)\right],\\
&&\\
B^{0,2}(u,v) &=&
\frac{1}{4}\left[
B(u,v) + J B(Ju,v) + J B(u,Jv) - B(Ju,Jv)
\right].
\end{eqnarray*}

In particular, note that $B^{1,1}(u,u)$ is the mean curvature vector of $\Sigma$ and so $B^{1,1}=0$ if and only if the surface is minimal. This decomposition of the second fundamental form is reminiscent of the splitting of the curvature tensor of a four-manifold as $\Rm = W^+ \oplus \Ric \oplus W^-$; the $B^{2,0}$ and $B^{0,2}$ components are actually conformally invariant, whilst $B^{1,1}$ depends on the metric itself.

Eells and Salamon observe that the lifted curve $\tilde f(\Sigma)$ is $J_+$-holomorphic if and only if $B^{0,2}=0$ and $J_-$-holomorphic if and only if $B^{1,1}=0$. (The condition $B^{2,0}=0$ corresponds to the lift being $J_+$-holomorphic for the opposite orientation on $X$.) Conversely, given an immersed curve $g\colon \Sigma \to Z$ which is everywhere transverse to the twistor fibration, composing with the projection $Z \to X$ gives an immersion $f \colon \Sigma \to Z$. If $g$ is $J_-$-holomorphic, say, then $g$ is the lift of $f$ and so $f$ is a minimal immersion. 

There is a standard way to adapt the above discussion to maps $f \colon \Sigma \to X$ which have certain isolated singularities.

\begin{definition}\label{branched}
A map $f \colon \Sigma\to X$ is called a \emph{branched immersion} if it is an immersion away from a finite set of points for each of which one can find a local coordinate in which $f$ has the form
$$
f(x,y) = \left(\re (x+iy)^N, \im (x+iy)^N, 0, 0\right) + O(x^N,y^N)
$$
in $C^1$. 
\end{definition}

One reason to make this definition is that these are precisely the singularities which occur when one considers minimal immersions (or immersions with $B^{0,2}=0$) with isolated singularities.

If $f$ is a branched immersion from an oriented surface, one can check that it is possible to define the tangent space to the image even at the branch points. This gives an oriented bundle $T_f\to\Sigma$, defined away from the branch points as before by $(T_f)_p = f_*(T_p\Sigma)$; moreover, $(T_f)_p$ varies continuously with $p$ even over the branch points. Taking the orthogonal complement in $TM$ gives a normal bundle $N_f$. Just as above, one can define the lift $\tilde f$ of $f$ to twistors space. The branch points correspond to points where the lift becomes tangent to a twistor line. 

This sets up the following analogue of the Weierstrass correspondence. (A ``non-vertical'' pseudoholomorphic curve is one which is not contained in a vertical twistor line.) 

\begin{theorem}[Eells--Salamon \cite{eells-salamon}]~
\begin{enumerate}
\item
There is a one-to-one correspondence between branched immersions of oriented surfaces in $X$ with $B^{1,1}=0$---i.e., branched minimal imm\-ersions---and non-vertical $J_-$-holomorphic curves in $Z$.
\item
There is a one-to-one correspondence between branched immersions of oriented surfaces in $X$ with $B^{0,2}=0$ and non-vertical $J_+$-holomorphic curves in $Z$.
\end{enumerate}
\end{theorem}

\subsection{An adjunction inequality}\label{adjunction}

When $J_+$ or $J_-$ is tamed by $\omega$, the fact that $\omega$ evaluates positively on pseudoholomorphic curves translates into a topological inequality for branched immersions whose lifts are pseudoholomorphic. 

Consider first the case of an embedding $f\colon \Sigma \to X$ of an oriented surface into a Riemannian four-manifold. As described above, the lift $\tilde f$ determines an almost complex structure on $f^*TX$ for which it splits as a sum of line bundles $T\Sigma \oplus N$. Upstairs in $Z$, this gives a splitting
$$
\tilde f^*TZ = T\Sigma \oplus N \oplus \tilde f^* V,
$$
where $V$ is the vertical tangent bundle of $Z$. We now evaluate the degree of $\tilde f^*TZ$ over $\Sigma$, recalling that $c_1(TZ)$ and $c_1(V)$ both depend on the choice of almost complex structure we use on $Z$. Just as in the calculation in \S\ref{cohomology}, for $J_+$, $c_1(TZ)=2c_1(V)$ (as originally observed by Hitchin \cite{hitchin2}) whilst for $J_-$, $c_1(TZ)=0$. Hence,
$$
\chi(\Sigma) + \Sigma\cdot \Sigma 
= 
\pm\int_\Sigma \tilde f^*c_1(V),
$$
with sign according to the choice of $J_\pm$ and where $\Sigma \cdot \Sigma$ denotes the self-intersection of $\Sigma$ in $X$. 

Assume now that the metric on $X$ satisfies the inequality (\ref{tame inequality}) of Theorem \ref{tame}  and $\det(W^++s/12)>0$, so that $\omega$ tames $J_+$ on $Z$. Then $c_1(V)=[\omega]$ and we deduce that whenever $f$ has $B^{0,2}=0$ (or, more generally, whenever $\tilde f$ is symplectic),
\begin{equation}\label{adjunction inequality for +ve case}
\chi(\Sigma) + \Sigma \cdot \Sigma > 0.
\end{equation}
On the other hand if $\det(W^++s/12) <0$, then $\omega$ tames $J_-$ and we deduce that whenever $f$ is minimal (or, again, whenever $\tilde f$ is symplectic),
\begin{equation}\label{adjunction inequality for -ve case}
\chi(\Sigma) + \Sigma \cdot \Sigma < 0.
\end{equation}

These adjunction inequalities should be seen as analogous to the fact that $\chi(\Sigma)<0$ for any minimally immersed surface in a negatively curved three-manifold. In this four-dimensional setting  inequalities such as (\ref{adjunction inequality for +ve case}) and (\ref{adjunction inequality for -ve case}) have appeared several times already in the literature. The archetypal example is that of a holomorphic curve in a complex surface with positive canonical bundle (in this context the inequality (\ref{adjunction inequality for -ve case}) follows from the adjunction formula; it is from here that the name ``adjunction inequality'' comes). In a Riemannian context, one of the earliest such results was due to Wang who showed that (\ref{adjunction inequality for -ve case}) held for embedded minimal surfaces in hyperbolic manifolds. This result was generalised by Chen--Tian \cite{chen-tian} who proved that (\ref{adjunction inequality for -ve case}) holds for embedded minimal surfaces in four-manifolds satisfying a variety of curvature inequalities. Several of these inequalities are subsumed and generalised by the curvature inequality (\ref{tame inequality}) in Theorem \ref{tame}. 

The above discussion also applies to branched immersions $f\colon \Sigma \to X$ of an oriented surface (Definition \ref{branched}). There is still a well defined tangent bundle $T_f$ and normal bundle $N_f$ to the image of $f$ and, just as above, we have
$$
c_1(T_f) + c_1(N_f)
=
\pm
\int_\Sigma \tilde f^*\, c_1(V),
$$
with sign according to the choice of $J_\pm$. The difference, however, is that $c_1(T_f)$ encodes both the topology of $\Sigma$ \emph{and} the branch points, whilst $c_1(N_f)$ encodes the topology of the embedding in a more complicated way than before. For an \emph{immersion}, for example, $c_1(N_f)$ includes both the self-intersection and the number of double points $d$, counted with sign:
$$
\int_\Sigma c_1(N_f) = f(\Sigma) \cdot f(\Sigma) - 2d.
$$

The effect of branch points is more complicated. In the complex case of a branched immersion of a holomorphic curve in a complex surface these phenomena are well understood. The Riemann--Hurwitz formula describes how branch points affect $c_1(T_f)$. Meanwhile, the change in $c_1(N_f)$ can be described via certain invariants of the link determined by each branch point. (The link of a branch point is obtained by considering a small ball in $X$ centred at the singularity and intersecting the boundary of the ball with the immersed surface.) These results have been generalised to the real case of branched immersions as in Definition \ref{branched}; the generalisation of the Riemann--Hurwitz formula is due to Gauduchon \cite{gauduchon2} whilst the  computation  of $c_1(N_f)$ in terms of link invariants is due to Ville \cite{ville2}.

\begin{definition}[Eells--Salamon \cite{eells-salamon}]
Given a branched immersion $f$ of an oriented surface into an oriented Riemannian four-manifold, the number $c_1(T_f)+c_1(N_f)$ is called the \emph{twistor degree} of $f$. (Note, Chen--Tian \cite{chen-tian} call this the \emph{adjunction number} by analogy with complex surfaces.)
\end{definition}

\begin{theorem}
Let $X$ be a compact Riemannian four-manifold whose curvature satisfies the inequality (\ref{tame inequality}) of Theorem \ref{tame}. Let $f\colon \Sigma \to X$ be a branched immersion. 
\begin{enumerate}
\item
If $\det(W^++s/12) >0$ and $f$ has $B^{0,2}=0$, then the twistor degree of $f$ is positive.
\item
If $\det(W^++s/12) <0$ and $f$ is minimal, then the twistor degree of $f$ is negative.
\end{enumerate}
\end{theorem}

\subsection{Gromov compactness}\label{gromov compactness}

The fundamental result in the theory of pseudoholomorphic curves in symplectic manifolds is Gromov's Compactness Theorem \cite{gromov}. Roughly speaking this states that if $(Z, \omega)$ is a compact symplectic manifold and $J$ an almost complex structure tamed by $\omega$ then any family of $J$-holomorphic curves $f\colon \Sigma \to Z$ for which the symplectic area $\int_\Sigma f^*\omega$ is independent of $f$ is weakly compact: any sequence $f_n$ of such curves has a subsequence which converges once one allows $\Sigma$ to develop ``bubbles'' (copies of $S^2$ glued to $\Sigma$ at those points where $\diff f_n$ tends to infinity) and nodes (corresponding to the degeneration of the conformal structure of $\Sigma$ induced by $f$). The limiting map $f_\infty$ has domain $\Sigma'$---the curve built from $\Sigma$ by adding bubbles and forming nodes---and $f\colon \Sigma'\to X$ is referred to as a ``cusp curve''. For brevity, we don't give a precise statement here and instead refer the reader to \cite{mcduff-salamon2} for a comprehensive discussion.

Combining Gromov compactness with the Eells--Salamon correspondence we see that for metrics satisfying the inequality (\ref{tame inequality}) of Theorem \ref{tame}, there is a similar compactness result for branched immersions which either have $B^{0,2}=0$ or are minimal, according to the sign of $\det(W^++s/12)$. To deduce the result stated below simply recall that for a branched immersion $f$, the symplectic area of the lift $\tilde f$ is, up to sign, exactly the twistor degree $c_1(T_f)+c_1(N_f)$.

\begin{theorem}\label{compactness theorem}
Let $X$ be a compact Riemannian four-manifold satisfying the curvature inequality (\ref{tame inequality}) of Theorem \ref{tame}. 
\begin{enumerate}
\item
If $\det(W^++s/12) >0$, then the set of branched immersions with $B^{0,2}=0$ and with fixed twistor degree is Gromov compact. I.e., every sequence of such maps has a subsequence which converges to a cusp curve.
\item
If $\det(W^++s/12) < 0$, then the set of branched minimal immersions with fixed twistor degree is Gromov compact.
\end{enumerate}
\end{theorem}

This should be seen as the analogue of the compactness result for minimal surfaces in a negatively curved three-manifold. What is perhaps a little surprising is that the curvature inequality (\ref{tame inequality}) is different from simply negative curvature. It is not hard to write down curvature tensors which satisfy (\ref{tame inequality}) with $\det(W^++s.12) < 0$ and yet have some sectional curvatures positive. 

As a contrast to Theorem \ref{compactness theorem} consider the example of the product of a hyperbolic surface $S$ with itself; this lies on the boundary of inequality (\ref{tame inequality}). (The 2-form $\omega$ is symplectic everywhere on $Z$ except along the equatorial circle separating the two sections of $Z\to M$ corresponding to plus and minus the K\"ahler form.) Take a sequence of increasingly long geodesics in $S$ which are all null-homologous. Their products give null-homologous minimal surfaces in $S \times S$; they  have twistor degree zero but area tending to infinity.

A important point to bear in mind is that the Eells--Salamon correspondence deals solely with \emph{non-vertical} curves. In the positive case, it may well happen that a limiting pseudoholomorphic curve has bubbles which are vertical or it may even be completely vertical itself. For example, in the case of the round metric on $S^4$, the twistor space is $\C\P^3$ with its standard K\"ahler structure. A generic linear $\C\P^1\subset \C\P^3$ (i.e., one transverse to the twistor fibration) will project to an embedded two-sphere in $S^4$ with $B^{0,2}=0$, but when one moves this $\C\P^1$ into a vertical twistor line, the projection in $S^4$ shrinks to a point. 

\section{Hyperbolic geometry and the conifold}\label{quadric cone}

Having considered mainly compact examples until now, we switch our focus to a particularly interesting non-compact symplectic manifold. As we have seen, the twistor spaces of $H^4$ and $H^2_{\C}$ carry natural symplectic structures. We begin in the first part of this section by showing that these are symplectomorphic. In fact, this symplectic manifold is already well-known; we show in the second part of this section that it is symplectomorphic to the total space of $\mathcal O(-1)\oplus \mathcal O(-1)\to\C\P^1$. The third part of this section explains how this picture gives ``one half'' of a hyperbolic interpretation of the threefold quadric cone $\{ xw-yz=0\}\subset \C^4$---also known as the conifold---and its desingularisations. 

\subsection{Paths of connections on $\SU(2)$}\label{paths of connections}

We begin by describing a framework in which it is easy to compare the definite connections coming from the symmetric spaces $S^4$, $\C\P^2$, $H^4$ and $H^2_\C$.

\subsubsection{A differential inequality}

Let $\epsilon_i$ be a standard basis for $\so(3)$ (i.e., with $[\epsilon_1, \epsilon_2] = \epsilon_3$ and cyclic permutations) and let $e_i$ be a standard basis for left-invariant 1-forms on $\SU(2)$ (i.e., with $\diff e_1 = e_2\wedge e_3$ and cyclic permutations). We consider a special class of connections on the trivial $\SO(3)$-bundle over $\SU(2)$, namely those of the form
$$
\nabla_B = \diff + \sum_{i=1}^3 a_i\,e_i\otimes \epsilon_i
$$
where the $a_i$ are constant. A path of such connections $B(r)$ defines a connection $\nabla_A = \frac{\del}{\del r} + \nabla_{B(r)}$ on $\R \times \SU(2)$. As we will see shortly, the definite connections coming from the symmetric spaces $S^4$, $\C\P^2$, $H^4$ and $H^2_\C$ all have this form.

Given a path $B(r)$, it makes sense to ask for the connection $A$ on $\R\times \SU(2)$ to be definite. This will translate into a differential inequality for the path $(a_1(r), a_2(r), a_3(r))$ in $\R^3$. If two paths can be deformed into each other whilst remaining definite, the corresponding symplectic forms are isotopic; it may then be possible to employ Moser's argument in this non-compact setting to show the manifolds are actually symplectomorphic. The first step is to compute the differential inequality for the path $(a_i)$.

Direct calculation gives that the curvature of $B$ is
$$
F_B
=
\begin{array}{r}
(a_1 + a_2a_3)\, e_2\wedge e_3 \otimes \epsilon_1\\
+(a_2 + a_3a_1)\, e_3\wedge e_1 \otimes \epsilon_2\\
+(a_3 + a_1a_2)\, e_3\wedge e_1 \otimes \epsilon_3
\end{array}
$$
(using $[\epsilon_1, \epsilon_2]=\epsilon_3$, $\diff e_1 = e_2\wedge e_3$ and cyclic permutations) whilst the curvature of $A$ is given by $F_A = \diff r \wedge\frac{\del B}{\del r} + F_B$ which is
$$
F_A
=
\begin{array}{r}
(a_1'\, \diff r \wedge e_1 + (a_1 + a_2a_3)\, e_2\wedge e_3)
\otimes \epsilon_1\\
+(a_2'\, \diff r \wedge e_2 + (a_2 + a_3a_1)\, e_3\wedge e_1)
\otimes \epsilon_2\\
+(a_3'\, \diff r \wedge e_3 + (a_3 + a_1a_2)\, e_1\wedge e_2)
\otimes \epsilon_3
\end{array}
$$
The following is immediate.
\begin{lemma}\label{definite path}
The path $(a_i)$ defines a definite connection if and only if the three quantities below are all non-zero and have the same sign:
$$
a_1'(a_1 +  a_2 a_3),\ 
a_2'(a_2 + a_3 a_1),\ 
a_3'(a_3 +  a_1 a_2).
$$
\end{lemma}

\subsubsection{Paths of left-invariant metrics}

We now describe how the four symmetric spaces $S^4$, $H^4$, $\C\P^2$ and $H^2_\C$ can be seen in this picture.  

We consider Riemannian metrics on $\R\times\SU(2)$ of the form
\begin{equation}\label{invariant metric}
\diff r^2 + \sum_{i=1}^3 f^2_i\, e_i^2
\end{equation}
where the $f_i$ are functions of $r$ and $\{e_i\}$ is a standard coframing for $\SU(2)$ as above. Such a metric determines connections on the $\SO(3)$-bundles $\Lambda^\pm$. On restriction to $\SU(2)\times\{r\}$, these bundles become naturally isomorphic to $T\SU(2)$. This is because of the standard fact that, given a non-zero $u\in \R^4$, $\Lambda^\pm\R^4$ is isomorphic to $\langle u\rangle^\perp$ via $v\mapsto (u\wedge v)^\pm$. Since we have chosen a (co)framing for $T\SU(2)$ we can see the connections on $\Lambda^\pm$ as paths of connections on the trivial bundle over $\SU(2)$. In fact, these connections are of exactly the form described above. 

\begin{lemma}
The Levi--Civita connection on $\Lambda^\pm$ of the metric (\ref{invariant metric}) is given by $\diff + \sum a_i e_i \otimes \epsilon_i$ where
\begin{eqnarray*}
a_1 & = &
\mp \frac{f_1'}{2}
+
\frac{f_1^2 - f_2^2 - f_3^2}{2f_2 f_3},\\
a_2 & = &
\mp \frac{f_2'}{2}
+
\frac{f_2^2 - f_3^2 - f_1^2}{2f_3 f_1},\\
a_3 & = &
\mp \frac{f_3'}{2}
+
\frac{f_3^2 - f_1^2 - f_2^2}{2f_1 f_2}.
\end{eqnarray*}
\end{lemma}
\begin{proof}
This is a direct calculation. It is made simpler by following a short-cut to compute the Levi--Civita connection on $\Lambda^+$, say, explained to us by Michael Singer. Let $\omega_i$ be a local orthonormal oriented basis for $\Lambda^+$ and write the connection in this basis as
$$
\nabla\left(
\begin{array}{c}
\omega_1\\
\omega_2\\
\omega_3
\end{array}
\right)
=
\left(
\begin{array}{ccc}
0 & + \alpha_3 \otimes \omega_2 & - \alpha_2 \otimes \omega_3\\
- \alpha_3\otimes \omega_1 & 0 & +\alpha_1\otimes \omega_3\\
+ \alpha_2 \otimes \omega_1 &  -\alpha_1 \otimes \omega_2 & 0
\end{array}
\right)
$$
for 1-forms $\alpha_i$. 

Since the $\omega_i$ are an orthonormal oriented basis, the maps $\alpha \mapsto \ast(\alpha \wedge \omega_i)$ define a triple of almost complex structures $J_i$ which satisfy the quaternionic relations $J_1J_2= J_3$ etc. Now the fact that the Levi--Civita connection is torsion free gives 
$$
\ast \diff \omega_1 = J_2(\alpha_3) - J_3(\alpha_2)
$$
along with similar equations for $\ast \diff \omega_2$ and $\ast \diff \omega_3$. Using the quaternionic relations one can solve for the $\alpha_i$, thus giving an explicit formula for the connection. E.g.,
$$
\alpha_1=\frac{1}{2}\left(
J_2(\ast \diff \omega_3) - J_3(\ast \diff \omega_2) - \ast \diff \omega_1 
\right).
$$

For the metric (\ref{invariant metric}) we do this calculation with the basis
\begin{eqnarray*}
\omega_1 &=&
f_1\,\diff r \wedge e_1 + f_2 f_3\,e_2 \wedge e_3,\\
\omega_2 &=&
f_2\,\diff r \wedge e_2 + f_3 f_1\,e_3 \wedge e_1,\\
\omega_3 &=&
f_3\,\diff r \wedge e_3 + f_1 f_2\,e_1 \wedge e_2.
\end{eqnarray*}

\end{proof}

For the round metric on $S^4$, we take all $f_i = \sin r$ for $r\in (0, \pi)$. Then the connection on $\Lambda^+$ is given by all $a_i= -(\cos r +1)/2$ whilst the connection on $\Lambda^-$ is given by all $a_i=(\cos r -1)/2$; by Lemma \ref{definite path} both give definite connections as expected. 

The Fubini--Study metric is given instead by $f_1=\sin r \cos r$ and $f_2 = f_3 = \sin r$ for $r \in (0, \pi/2)$. The connection on $\Lambda^+$ has $a_1= -(\cos^2 r +1)/2$ and $a_2=a_3=-\cos r$ which is definite as expected. The connection on $\Lambda^-$ meanwhile has $a_2=a_3=0$ and so is not definite. (We have implicitly given $\C\P^2$ the non-standard orientation in this calculation.)

The hyperbolic metric is given by all $f_i = \sinh r$ for $r \in (0,\infty)$. The connection on $\Lambda^+$ has all $a_i =-(\cosh r +1)/2$ whilst that on $\Lambda^-$ has all $a_i=(\cosh r -1)/2$. Again we confirm that both are definite. 

Finally, the complex-hyperbolic metric is given by $f_1=\sinh r \cosh r$ and $f_2 = f_3 = \sinh r$ for $r\in (0,\infty)$. The connection on $\Lambda^+$ has $a_1 =  - (\cosh^2 r +1)/2$ and $a_2=a_3=- \cosh r$, which is definite. The connection on $\Lambda^+$ has $a_2=a_3=0$ and so is not definite. (Again, $H^2_\C$ has the non-standard orientation here.)

With this description of the definite connections in hand we can prove the following.

\begin{proposition}
The twistor spaces of $H^4$ and $H^2_\C$ with their natural symplectic structures are symplectomorphic.
\end{proposition}

\begin{proof}
We first show that the two symplectic forms are isotopic by considering the linear isotopy between their corresponding paths. This is given by
\begin{eqnarray*}
a_1(r,t) &=&
-\frac{1}{2}(1 + \cosh r) + \frac{t}{2}\cosh r(1 - \cosh r),\\
a_2(r,t) = a_3(r,t) &=&
-\frac{1}{2}(1 + \cosh r) + \frac{t}{2}(1 - \cosh r).
\end{eqnarray*}
The path $a_i(r,0)$ corresponds to $H^4$ whilst the path $a_i(r,1)$ corresponds to $H^2_\C$. It is straightforward to check from Lemma \ref{definite path} that for each $t\in [0,1]$, the path $a(r,t)$ gives a definite connection and hence a symplectic form $\omega_t$ interpolating between the symplectic forms corresponding to $H^4$ and $H^2_\C$.

To turn the isotopy into a symplectomorphism we use Moser's argument. The first step is to write $\omega_t = \omega_0 + \diff \alpha_t$ for a path of 1-forms $\alpha_t$ on the sphere bundle $Z$. Now,  $\omega_t$ is the curvature of the connection $\nabla_t$ on the vertical tangent bundle $V$ determined by the $\SO(3)$-connection $A_t$:
$$
A_t = \diff + \sum a_i(r,t) e_i \otimes \epsilon_i
=
A_0 + tC
$$
where 
$$
C(r)=\frac{1}{2}(1 - \cosh r)(\cosh r\, e_1\otimes \epsilon_1 + e_2\otimes \epsilon_2 + e_3 \otimes \epsilon_3)
$$
is independent of $t$. Recall that $\so(E)$-valued $p$-forms on the four-manifold determine genuine $p$-forms on $Z$ (via pulling back and pairing with the fibrewise moment map $Z^2\to \so(E)^*$). We denote by $c\in \Omega^1(Z)$ the 1-form corresponding in this way to $C$. It is straight forward to verify that upstairs on $Z$ the connections in $V$ are related by $\nabla_t = \nabla_0 + tc$ and hence $\omega_t = \omega_0 + t \diff c$. 

Using the symplectic form $\omega_t$ to identify $T^*Z \cong TZ$ gives a family of vector fields $u_t$ solving $\iota_{u_t}\omega_t = c$. To complete the argument we must show that the flow of $u_t$ exists for $t\in[0,1]$. In other words, given any $p\in Z$, we claim that the solution to the ODE $\gamma'(t) = u_t(\gamma(t))$ with $\gamma(0)=p$, defined initially for $t\in[0, T)$ with $T \leq 1$, travels a finite distance. Completeness will then ensure the flow can be continued past $t=T$ and up to $t=1$. To prove this claim we will control the radial component of the image $\pi_*(u_t)$ of the vector field downstairs.

The 1-form $c$ is horizontal, hence we only need consider the horizontal components of $\omega_t$---i.e., the curvature of $A_t$---when computing $u_t$. When approximating the behvaiour of tensors for large values of $r$, we will write $P \sim Q$ to mean that the coefficients of $P$ (with respect to the bases built from $\diff r, e_i, \epsilon_i$) are controlled above and below, independently of $r$, by constant multiples of the coefficients of $Q$. So $C\sim \hat C$ where
$$
\hat 
C
= 
e^{2r}e_1\otimes \epsilon_1
+e^{r} e_2\otimes \epsilon_2
+e^{r} e_3\otimes \epsilon_3,
$$
whilst, for any $t>0$, $F_{A_t} \sim \hat F$ where
$$
\hat F
=
\begin{array}{r}
\left(
e^{2r}\,\diff r\wedge e_1 + e^{2r}\, e_2\wedge e_3
\right)\otimes\epsilon_1\\
+
\left(
e^r\,\diff r\wedge e_2 + e^{3r}\, e_3\wedge e_1
\right)\otimes\epsilon_2\\
+
\left(
e^r\,\diff r\wedge e_3 + e^{3r}\, e_1\wedge e_2
\right)\otimes\epsilon_3
\end{array}
$$
where the control is uniform in $t$ also, provided $t>t_0>0$ is bounded strictly away from zero. (The control depends on $t_0$ and, in fact, deteriorates as $t_0 \to 0$, but this will not affect us.)

Now $\hat C$ determines a 1-form $\hat c$ on $Z$, just as $C$ determines $c$; also, $\hat F$ determines a 2-form $\hat \omega$ on $Z$ in the same way that $F_{A_t}$ determines the purely horizontal component of $\omega_t$. Let $\xi$ denote any lift of $\del_r$ to $Z$ (i.e., any vector such that $\pi_*(\xi)=\del_r$. Note that $\iota_{\xi} \hat \omega = \hat c$. It follows that, provided $t$ is bounded strictly away from zero, the radial component of $\pi_*(u_t)$ is bounded independently of $r$ and $t$. 

To complete the proof, assume the flow starting at some point $\gamma(0)$ exists for $t\in [0,T)$. Since the radial component of $\pi_*(u_t)$ is controlled uniformly in $t$ and $r$ for $t>T/2$ we see that the radial distance travelled, as $t$ runs from $T/2$ to $T$ is finite.
\end{proof}

\subsection{Symplectic geometry of the small resolutions}\label{symplecto to small res}

We now explain how the twistor space of $H^4$ is a familiar symplectic manifold, namely the total space of $\mathcal O(-1)\oplus \mathcal O(-1)\to \C\P^1$. It will be convenient to view this bundle as the small resolution of the conifold $xw-yz=0$ in $\C^4$ in the manner described below. 

\subsubsection{The small resolutions}

The conifold and its desingularisations have been much studied by both mathematicians and physicists. We learnt about this subject from \cite{smith-thomas-yau} and refer the interested reader to that article and the references therein. 

Let $Q = \{xw- yz=0\}$ denote the conifold in $\C^4$. To resolve the double-point $0\in Q$ begin by blowing up the origin in $\C^4$. The proper transform $\hat Q$ of $Q$ intersects the exceptional $\C\P^3$ in the quadric surface given by the same equation as $Q$. The maps $[x:y:z:w] \mapsto [x:y]=[z:w]$ and $[x:y:z:w]\mapsto [x:z] = [y:w]$ show that the quadric surface is biholomorphic to $\C\P^1\times \C\P^1$; blowing down one of the rulings inside of $\hat Q$ gives a resolution $p_+\colon R^+\to Q$. Away from $0\in Q$, $p_+$ is an isomorphism; the fibre over $0$  is a copy of $\C\P^1$ which we denote $\C\P^1_+$. Blowing down the other ruling results in a resolution $p_-\colon R^-\to Q$ with exceptional curve $\C\P^1_-$. 

The map $(x,y,z,w)\mapsto [x:y]=[z:w]$, defined initially on $Q \setminus 0$, extends to one resolution whilst $(x,y,z,w)\mapsto [x:z] = [y:w]$ extends to the other, giving projections $q_\pm\colon R^\pm\to\C\P_\pm^1$. These show that each $R^\pm$ is abstractly biholomorphic to the total space of $\mathcal O(-1)\oplus \mathcal O(-1)$ but it is important to note that they are not identifiable in a way which respects the projections $p_\pm$. 

The group $\SO(4,\C)$, of linear isomorphisms of $\C^4$ preserving the quad\-ratic form $xw-yz$, acts on $Q$ and the action lifts to the two resolutions $R^\pm$. One way to describe these lifts is to view $\C^4$ as $\End \C^2$ and the quadratic form as the determinant. Pre- and post-multiplication by elements of $\SL(2,\C)$ preserve the quadratic form and this sets up an isomorphism
$$
\SO(4,\C) \cong \frac{\SL(2,\C) \times \SL(2,\C)}{\pm1}
$$
In the action of $\SO(4,\C)$ on $R^+\cong \mathcal O(-1) \oplus \mathcal O(-1)$, one copy of $\SL(2,\C)$ acts via the natural action of $\SL(2,\C)$ on the total space of $\mathcal O(-1)$. The other copy acts linearly on each fibre of $\mathcal O(-1)\oplus \mathcal O(-1)$ covering the identity on $\C\P^1$; it does this via the natural trivialisation $\End(\mathcal O(-1)\oplus \mathcal O(-1))\cong \End \underline{\C^2}$. In the action on $R^-$ the r\^{o}les of the two factors are swapped.

The small resolutions are, of course, K\"ahler; they carry a one-parameter family of compatible symplectic forms described as follows. The pull-back of the standard symplectic structure on $\C^4$ via $p_\pm\colon R^\pm \to Q\subset \C^4$ is degenerate on the exceptional $\C\P^1_\pm$. To remedy this, use the projection map $q_\pm\colon R^\pm\to \C\P_\pm^1$ to pull back the Fubini--Study form. This gives a symplectic form on $R^\pm$; define
$$
\omega = p_\pm^*\omega_{\C^4} + q_\pm^*\omega_{\mathrm {FS}}.
$$

\subsubsection{The small resolutions as twistor spaces.}

We now describe the small resolutions as the two twistor spaces of $\R^4$. 

Take $\R^4$ with its standard Euclidean structure and write $\C^4 \cong \R^4 \otimes \C$; now take the complex quadratic form to be the complexification of the Euclidean quadratic form on $\R^4$. I.e., we have chosen coordinates in $\C^4$ so that $Q$ is given by $\{\sum z_j^2 = 0\}$. Note that $\SO(4)$ acts on $\C^4$, extending the usual action on $\R^4$ by complex linearity and giving an embedding $\SO(4)\subset \SO(4,\C)$. This action lifts to the small resolutions.

\begin{lemma}
The group $\SO(4)$ acts by symplectomorphisms on the small resolutions.
\end{lemma}
\begin{proof}
The action of $\SO(4)$ on $\C^4$ is unitary; indeed it commutes with conjugation and so preserves the Hermitian inner-product $\omega_{\C^4}$ which is built by combining conjugation and the complex bilinear form (in fact $\SO(4) = U(4)\cap\SO(4,\C)$). Hence the lift of the action of $\SO(4)$ to the small resolutions preserves $p_\pm^*\omega_{C^4}$.

The action on the small resolutions preserves the exceptional curves. The $\C\P^1_\pm$ inherit their Fubini--Study metrics from the unitary structure on $\C^4$. Since $\SO(4)$ acts unitarily on $\C^4$, the induced actions on $\C\P^1_\pm$ are also isometric. Now, as $q_\pm$ is equivariant, it follows that the action of $\SO(4)$ on the small resolutions also preserves $q^*_\pm\omega_{\mathrm{FS}}$ and hence $\omega$.
\end{proof}

\begin{lemma}
There is an $\SO(4)$-equivariant identification of $R^\pm$ with the unit sphere bundle in $\Lambda^\pm\to\R^4$. 
\end{lemma}
\begin{proof}
To begin with, consider the isometric action of $\SO(4)$ on the exceptional $\C\P^1_\pm$. This gives a map $\SO(4)\to\SO(3)_+\times \SO(3)_-$; but there is already a standard such map induced by the action of $\SO(4)$ on $\Lambda^2\R^4=\Lambda^+\R^4\oplus \Lambda^-\R^4$. So we can identify $\C\P^1_\pm$ with the unit spheres in $\Lambda^\pm$ in an $\SO(4)$-equivariant fashion.

To see the rest of the twistor fibration, write out the real and imaginary parts of the defining equation for $Q$:
$$
Q = \{ u+iv : |u| = |v|,\ \langle u, v\rangle = 0\}
$$
where $u, v \in \R^4$ and $|\cdot |$ and $\langle \cdot, \cdot \rangle$ are the Euclidean structures on $\R^4$. Let $t \colon Q\to \R^4$ be given by taking the real part, $t(z)=\re (z)$. At $u\neq 0 \in\R^4$, the fibre of $t$ is the sphere of radius $|u|$ in the three-plane $\langle u \rangle ^\perp$. This sphere can be identified with the unit sphere in $\Lambda^+\R^4$ by $v \mapsto \frac{(u\wedge v)^+}{|(u\wedge v)^+|}$. Away from the origin this identifies $t\colon Q\setminus 0 \to \R^4\setminus 0$ with the twistor projection in an $\SO(4)$-equivariant manner. Moreover, on restriction to $Q\setminus 0$, the projection $q_+$ identifies each fibre of $t$ with $\C\P^1_+$. Hence $t$ extends to a map $t_+\colon R^+\to \R^4$ completing the identification with the unit sphere bundle in $\Lambda^+\to \R^4$.

Of course, we could have used the anti-self-dual part of $u\wedge v$ to identify $Q\setminus 0$ with the opposite twistor space of $\R^4 \setminus 0$ (i.e., the unit sphere bundle in $\Lambda^-$) and this gives an identification of the other small resolution with the other twistor space.
\end{proof}

One might ask how this tallies with the fact that the twistor spaces are biholomorphic to $\mathcal O(1) \oplus \mathcal O(1)$. The answer, of course, is that the diffeomorphisms described above are not biholomorphisms. The vector bundle projections $q_\pm$ are transverse to $t_\pm$ and the fibres of $q_\pm$ define sections of the twistor spaces which correspond to linear complex structures on $\R^4$ compatible with the Euclidean structure. However, the fibres of $q_\pm$ already have linear complex structures by virtue of mapping isomorphically to affine subspaces under $p_\pm\colon R^\pm \to Q\subset\C^4$ and these complex structures are \emph{minus} those arising in the twistorial picture. 

In summary, the small resolutions admit three projections: $p_\pm \colon R^\pm \to Q$; $q_\pm \colon R^\pm \to \C\P^1_\pm$ and $t_\pm \colon R^\pm \to \R^4$. $\SO(4)$ acts on all these spaces making all the projections equivariant. The fibres of $q_\pm$ and $t_\pm$ are transverse and, if we think of the fibres of $t\pm$ as the unit spheres in $\Lambda^\pm$, $q_\pm$ maps each of them isometrically onto $\C\P^1_\pm$. Conjugation on $Q$ lifts to an involution on $R^\pm$ which is the antipodal map on the fibres of $t_\pm$---the usual ``real structure'' of twistor theory.

\subsubsection{The small resolution via a definite connection}

Symplectically then, the small resolution, $R^+$ say, can be thought of as an $\SO(4)$-invariant symplectic form on the twistor space $Z$ of $\R^4$. Moreover, under the identification $R^+ \cong Z$, the antipodal map on the fibres of $Z$---the usual ``real structure'' of twistor theory---corresponds to the action of conjugation on $Q$ lifted to $R^+$. This swaps the sign of both $\omega_{\mathrm{FS}}$ and $\omega_{\C^4}$ and hence also of $\omega$. We will now show that this, combined with the fact that $(R^+, \omega)$ has infinite symplectic volume, is enough to determine $\omega$ completely.

\begin{proposition}
Let $\omega$ be an $\SO(4)$-invariant symplectic form with infinite total volume on the twistor space $t\colon Z\to\R^4$. Suppose, moreover, that $\omega$ changes sign under the fibrewise antipodal map and gives the fibres area $4\pi$. Then there is an $\SO(4)$-equivariant identification $H^4\to\R^4$ such that $\omega$ is the symplectic form associated to the hyperbolic metric.
\end{proposition}

\begin{proof}
Let $u \in \R^4\setminus 0$ and let $\SO(3) \subset \SO(4)$ denote the stabiliser of $u$. The restriction of $\omega$ to the fibre $t^{-1}(u)$ is an $\SO(3)$-invariant 2-form on $S^2$ and so is, up to a sign, the standard area form. As the fibrewise restriction of $\omega$ is nondegenerate, the symplectic complements to the fibres define a horizontal distribution $H$ in $Z$. We will show that $H$ is the definite connection coming from the hyperbolic metric.

Let $H_\infty$ be the integral horizontal distribution in $Z\to \R^4$ determined by the flat Levi--Civita connection of the Euclidean metric on $\R^4$. The difference $H - H_\infty$ of the two connections is a 1-form $a$ on $\R^4$ with values in vertical vector fields. As the $\SO(4)$-action preserves both $H$ and $H_\infty$, $a$ is $\SO(4)$-equivariant. Note that even though the connection $H_\infty$ preserves the natural metrics on the fibres of $Z\to \R^4$, we don't yet know that $H$ does. This will follow from an explicit description of $a$. 

Given a point $u\in \R^4 \setminus 0$, we interpret $a$ as a linear map
$$
a \colon T_u\R^4 \to C^\infty (TS^2_u),
$$
where $S^2_u$ is the unit sphere in the three-plane orthogonal to $u$. The key step is the following claim:

\emph{Claim. The 1-form $a$ vanishes on vectors parallel to $u$. On the three-plane orthogonal to $u$, $a(v)$ is a scalar multiple of the unit speed rotation of $S^2_u$ about $v$. Moreover, the scalar factor does not depend on $v$ (although it does depend on $|u|$).}

To prove the claim, let $\SO(3)\subset\SO(4)$ denote the stabiliser of $u$. We will exploit the $\SO(3)$-equivariance of $a\colon T_u\R^4\to C^\infty(TS^2_u)$. First consider $a(v)$ when $v$ is parallel to $u$; the $\SO(3)$-action fixes $u$ and $v$ and so must also fix $a(v)$. However, the only $\SO(3)$-invariant vector field on $S^2$ is zero, hence $a(v)=0$. 

Next we consider the case when $v$ is perpendicular to $u$. Restricting $a$ to the three-plane perpendicular to $u$ gives a linear map $a\colon \langle u \rangle^\perp \to C^\infty(TS^2_u)$. By identifying each $v\in\langle u \rangle^\perp$ with unit speed rotation of $\langle u \rangle^\perp$ about $v$, we can consider $a$ as an $\SO(3)$-equivariant map $\so(3) \to C^\infty(TS^2)$. In fact, $a$ is also invariant under the antipodal map and so is $\mathrm{O}(3)$-equivariant. The claim now follows from the classical fact that the only $\mathrm{O}(3)$-invariant maps $\so(3)\to C^\infty(TS^2)$ are multiples of the standard embedding.

Since $H_\infty$ preserves the metrics on the fibres of $t\colon Z\to \R^4$ it follows from the italicized claim that $H$ does too; in other words the horizontal distribution of $\omega$ defines an $\SO(3)$-connection. This connection in turn defines a 2-form on the sphere bundle $Z$, but it is more-or-less tautological that this 2-form is just $\omega$ again, hence $\omega$ is induced by the definite connection $H$.

Next, consider the restriction of $H$ to the sphere $S^3$ of radius $r$ in $\R^4$. At each point $u\in S^3$, the fibre of $t$ is identified with the unit sphere in $T_uS^3$ and so on restriction $Z|_{S^3}$ is the unit sphere bundle in $TS^3$. To trivialise this bundle, fix an $\SU(2)\subset\SO(4)$ and let $\epsilon_i$ be an oriented orthonormal basis for $\su(2)$. Their images on $S^3$ trivialise $TS^3$ and so also $Z|_{S^3}$; the claim proved above says that in this trivialisation, the restriction of $H$ is given by 
$$
\nabla_{H_{r}} = \diff + C(r)\sum_{i=1}^3 e_i\otimes \epsilon_i
$$
for some constant $C(r)$, where $e_i$ is the framing of $T^*S^3$ dual to that induced by the $\epsilon_i$.

Moreover, the first part of the claim says that $H$ is in ``radial gauge''; that is to say it is completely determined by its restriction to the concentric spheres:
$$
\nabla_{H} = \frac{\del}{\del r} + \nabla_{H_{r}}
$$
(This is often called ``temporal gauge'' for connections on a cylinder $Y\times \R$.) In other words, we have written $H$ as a path of connections in precisely the form considered in \S\ref{paths of connections} with all $a_i = C(r)$. 

It remains to describe $C(r)$. First, since the connection extends over the origin in $\R^4$, $C(0)=0$ or $C(0)=-1$ (recall from the picture in \S\ref{paths of connections} that these are the only choices which give a flat connection on the trivial bundle over $S^3$). Secondly, the connection is definite and has infinite symplectic volume. It follows from the description in \S\ref{paths of connections} (in particular, Lemma \ref{definite path}) that $C(r)$ is strictly monotonic in $r$ and either gives a diffeomorphism $[0, \infty) \to [0,\infty)$ or $[0,\infty)\to [-1,-\infty)$. In either case, after appropriately rescaling the radial coordinate, this is the same as one of the definite connections coming from the hyperbolic metric.
\end{proof}

\subsection{Hyperbolic geometry and the model conifold transition}\label{mirror symmetry}

The final part of this section is of a highly speculative nature. We hope that despite this, it is worth writing down.

The preceding discussion shows that isometries of hyperbolic four-space act in a canonical way by symplectomorphisms on the small resolutions of the conifold. There is another way to desingularise the conifold, namely by smoothing it to the variety $S = \{xw-yz=1 \}$. As we will explain, isometries of hyperbolic \emph{three}-space (almost) act by \emph{biholomorphisms} on the smoothings $S$. These two different ways to desingularise $Q$ have been the subject of much work by both mathematicians and physicists. The passage from smoothing to resolution is called a ``conifold transition'' in the physics literature. As is described in various articles in \cite{essays}, the smoothing and small resolution are considered by physicists to be a ``mirror pair''. Since mirror symmetry swaps symplectic and complex, it is perhaps not too ridiculous to hope that there may be some duality between three- and four-dimensional hyperbolic geometry. In any case, there are certainly similarities between the two situations and these alone make it useful to consider both pictures side-by-side. 

In the complex picture, a natural starting point is to look for compact complex quotients of $S$. It is here that hyperbolic three-manifolds provide an answer. Note that, viewing the complex quadratic form on $\C^4\cong \End\C^2$ as the determinant, $S =\SL(2,\C)$. Meanwhile, $H^3$ can be seen as the symmetric space $\SL(2,\C)/\SU(2)$ and the resulting $\SU(2)$-bundle $\SL(2,\C)\to H^3$ is just the spin bundle of $H^3$. Now $\PSL(2,\C)$ is the isometry group of $H^3$ and so a compact hyperbolic three-manifold is determined by a cocompact lattice $\Gamma \subset\PSL(2,\C)$. It is a standard fact that this lattice lifts to an inclusion $\iota \colon \Gamma \to \SL(2,\C)$ and we can take the quotient $\SL(2,\C)/\Gamma$. This is a compact complex quotient of $S$ which is simply the spin bundle of the original hyperbolic three-manifold. Superficially at least, there is some similarity between these complex threefolds and the twistor spaces of hyperbolic four-manifolds: both are sphere bundles over hyperbolic manifolds naturally derived from spinor geometry. (We are grateful to Maxim Kontsevich for making us aware of these examples and also for directing us to the article of Ghys \cite{ghys}.)

These are not the only compact quotients of $S$, however. Whilst Mostow rigidity ensures that the lattice $\Gamma$ is rigid in $\SL(2,\C)$ it is certainly not rigid in $\SO(4,\C)$ (at least when $H^1(M)\neq 0$). Ghys shows in \cite{ghys} that if $\rho \colon \Gamma \to \SL(2,\C)$ is any homomorphism sufficiently close to the trivial homomorphism then the embedding $\Gamma \to \SO(4,\C)$ given by
$$
\gamma \mapsto [\rho(\gamma), \iota(\gamma)]
$$
leads to a different smooth complex quotient $S/\Gamma$ which is diffeomorphic to the spin bundle of the hyperbolic three-manifold but with a genuinely different complex structure. In a forthcoming paper \cite{fine-panov} we will show that in fact this result holds for any $\rho$ lying in an open set containing all homomorphisms with image in $\SU(2)$. Moreover, all the quotients have holomorphically trivial canonical bundle, making them ``holomorphic Calabi--Yaus''. Intriguingly, the inequality defining this open set of homomorphisms is founded on exactly the same linear algebra underlying the definition of a definite connection, namely that of maximal definite subspaces in $\Lambda^2$. 

It is natural to wonder if there is any deeper relationship between these symplectic and complex pictures, besides the immediate superficial similarity outlined here.

\section{Concluding remarks}\label{concluding remarks}

Even leaving aside the speculation of the previous section, it seems, to us at least, that there are many concrete unanswered questions concerning definite connections. Perhaps the most important is to determine whether or not the only closed four-manifolds which admit \emph{positive} definite connections are $S^4$ and $\C\P^2$. In fact, either outcome here would be interesting. On the one hand, if $S^4$ and $\C\P^2$ are the only possible bases for a positive definite connection this would amount to a gauge theoretic ``sphere theorem''. There are many results in Riemannian geometry going back several years which ensure that only the sphere admits a metric satisfying various curvature inequalities. If one could show that only $S^4$ and $\C\P^2$ admit positive definite connections, not only would this imply a new Riemannian ``sphere'' theorem in the classical sense (only $S^4$ and $\C\P^2$ admit metrics with $\mathcal D>0$ and $\det(W^++s/12)>0$), it would also, to our knowledge at least, be the first such result phrased purely in terms of a connection on some auxiliary bundle.

On the other hand, if there is another closed four-manifold admitting a positive definite connection this would give a symplectic Fano manifold. Unlike for $S^4$ and $\C\P^2$, however, this Fano would not be K\"ahler, making it the first such symplectic manifold known.  To see that it would not be K\"ahler, it suffices to look at the classification of Fano K\"ahler threefolds. Fanos arising from definite connections have even index; in the K\"ahler case, this leads to $\C\P^3$ or a del Pezzo threefold of index equal to 2. These have been classified by Iskovskih (see, for example, Theorem 3.3.1 of \cite{iskovskih}). The formulae in \S\ref{cohomology} imply that the bundle $L$ with $L^2=-K_Z$ satisfies $L^3 = 2(2\chi +3 \tau)$. There are five varieties on Iskovskih's list for which $L^3$ is even. Those with $L^3=2$ and $L^3=4$ have $b_2(Z)=1$, meaning the underlying four-manifold would have $b_2(M)=0$ and so $2\chi+3\tau=4$ contradicting the value of $L^3$. One of the candidates with $L^3=6$ is $\C\P^1\times \C\P^1\times \C\P^1$; the underlying four-manifold would then have $\chi = 4$ and so $b_2(M)=2$, $\tau$ equal to one of $2, -2,0$ all of which contradict $L^3=6$. The other del Pezzo is the flag manifold.

There are also questions concerning the negative case. Foremost amongst these is to understand whether or not a compact simply connected example exists. If not, what constraints does the existence of a negative definite connection place on the fundamental group? In a similar vein, all the compact examples presented here come via metrics with negative sectional curvatures and hence have vanishing higher homotopy groups. Are there any examples with non-vanishing higher homotopy? 

One concrete way to approach this problem could be to desingularise isolated orbifold points. For example, consider the action of $\Z_n$ on $\C^2$ generated by multiplication by an $n^\mathrm{th}$-root of unity. The minimal resolution $\mathcal O(-n) \to \C^2/\Z_n$ contains an essential sphere; moreover, identifying $\C^2 \cong H^4$ in a way that the action preserves the hyperbolic metric, we see that the quotient carries a hyperbolic orbifold metric. The question is, does the resolution carry a metric with $\mathcal D>0$ and $\det(W^++s/12)$? 

In this non-compact model situation, it is straightforward to see that such a metric exists. In the notation of \S\ref{paths of connections}, consider the metric on $\C^2\setminus 0 = \R_{>0} \times \SU(2)$ given by
$$
g_n = \diff r^2 + n\sinh^2\! r\, e_1^2 + n\cosh^2\! r\,(e_2^2 + e_3^2),
$$
where the $e_i$ are a standard basis for $\su(2)^*$. This metric has the correct singularity at $0$ to ensure that the pull-back to $\mathcal O(-n)$ extends smoothly over the exceptional curve.
As in \S\ref{paths of connections} the Levi--Civita connections on $\Lambda^\pm$ can be explicitly computed. Lemma \ref{definite path} tells us that the connection on $\Lambda^+$ is never definite; however, for $n\geq 3$, the connection on $\Lambda^-$ is negative definite. This example shows that, in the non-compact setting at least, there is a distinct difference between metrics with $\mathcal D>0$ and $\det(W^++s/12)<0$ and metrics of negative sectional curvature.

Incidentally, the inequality $n\geq 3$ is actually an instance of the adjunction inequality (\ref{adjunction inequality for -ve case}) of \S\ref{adjunction}. Any metric on $\mathcal O(-n)$ pulled back from an $\SU(2)$-invariant metric on $\C^2\setminus 0$ makes the exceptional curve totally geodesic and so, in particular, minimal. If the metric also satisfies the inequality (\ref{tame inequality}) of Theorem \ref{tame}, with $\det(W^++s/12)<0$, then the adjunction inequality applied to the exceptional curve gives $n \geq3$.

To produce compact examples, one might attempt to use this model to resolve the orbifold points of a compact hyperbolic orbifold with the correct type of isolated singularities. It is here that the fact that $g_n$ only has $\mathcal D>0$ in one orientation causes problems. For example, suppose one had a $\Z_n$-action on a hyperbolic four-manifold with isolated fixed points; suppose, morever, that at each fixed point $x$ it were possible to identify $T_xX \cong\C^2$ in such a way that the action is given by multiplication by an $n^\mathrm{th}$ root of unity. One could then resolve all the singularities by gluing in copies of $g_n$; however in order to produce a metric with $\mathcal D>0$, the distinguished orientations induced at each fixed point via $T_xX \cong \C^2$ must be coherent. In other words, there must be a global choice of orientation which agrees at each fixed point with the distinguished orientation there. 

One manifold to which it may be possible to apply this construction is the Davis manifold. (There is a detailed description of this hyperbolic manifold and its isometries in \cite{ratcliffe-tschantz}.) It is possible to write down various cyclic group actions; unfortunately, for the actions we found, the induced orientations at the fixed points were never coherent. Incidentally, attempting to understand this phenomenon revealed to us the paucity of concrete examples of hyperbolic four-manifolds. 

Returning to the general question of existence of definite connections, the only obstruction we know is $2\chi + 3\tau >0$, but presumably there are others and the class of manifolds admitting, say, negative definite connections is limited in some way. For example, if $S^4$ admitted a negative definite connection this would give $\C\P^3$ the structure of a symplectic Calabi--Yau, which seems counter-intuitive.

In the purely Riemannian setting there are further questions one could consider. For example, the space of algebraic curvature tensors for which $\mathcal D$ is definite has six components and we have only exhibited complete examples in two of them, namely those for which $\mathcal D >0$ and $W^+ + s/12$ is definite. Do there exist examples with $\mathcal D<0$ or with $W^++s/12$ indefinite?

We hope that this discussion, along with the various questions and suggestions in the text, shows that there is ample scope for future work in this area.

\bibliographystyle{plain}
\bibliography{symp_CYs_bib}

{\small \noindent {\tt joel.fine@ulb.ac.be }} \newline
D\'epartment de Math\'ematique,
Universit\'e Libre de Bruxelles CP218,\\ 
Boulevard du Triomphe,
Bruxelles 1050,
Belgique.\\
{\small \noindent {\tt d.panov@imperial.ac.uk}} \newline
Department of Mathematics,
Imperial College,\\
London SW7 2AZ,
UK.

\end{document}